\documentclass[twoside]{article}
\usepackage{amssymb}
\usepackage{indentfirst}
\usepackage{dsfont}
\usepackage{ntheorem}
\usepackage{amsmath}

\begin{document}
\baselineskip 13.5pt \noindent \thispagestyle{empty}

\markboth{\centerline{\rm Z. Tu \; \& \;L. Wang }}{\centerline{\rm
Proper mappings between Hua domains}}

\begin{center} \Large{\bf  Rigidity of proper holomorphic mappings between}\end{center}
\begin{center} \Large{\bf equidimensional Hua domains}\end{center}

\vskip 8pt

\begin{center}
\noindent\text{Zhenhan TU\; and \; Lei WANG$^{*}$}\\
\vskip 8pt \noindent\small {School of Mathematics and Statistics,
Wuhan University, Wuhan, Hubei 430072, P.R. China} \\
\noindent\text{Email: zhhtu.math@whu.edu.cn (Z. Tu),\;}
{wanglei2012@whu.edu.cn (L. Wang)}
\renewcommand{\thefootnote}{{}}
\footnote{\hskip -16pt {$^{*}$Corresponding author. \\ } }
\end{center}

\vskip 8pt

\begin{center}
\begin{minipage}{13cm}
{\bf Abstract.}
 {\small Hua domain, named after Chinese mathematician Loo-Keng Hua,
 is defined as a domain in $\mathbb{C}^{n}$ fibered over an irreducible bounded symmetric
 domain $\Omega\subset \mathbb{C}^{d}\;(d<n)$ with the fiber over
 $z\in \Omega$ being a $(n-d)$-dimensional generalized complex ellipsoid $\Sigma(z)$.
In general, a Hua domain is a nonhomogeneous domain without smooth
boundary. The purpose of this paper is twofold. Firstly, we obtain
what seems to be the first rigidity results on proper holomorphic
mappings between two equidimensional Hua domains.
 Secondly, we determine the explicit form of the biholomorphisms
 between two equidimensional Hua domains. As a special conclusion  of this paper, we completely describe the
group of holomorphic automorphisms of the Hua domain.

}
\end{minipage}
\end{center}

\newtheorem{lem}{Lemma}[section]
\newtheorem{thm}{Theorem}[section]
\newtheorem{prop}{Proposition}[section]

\section{Introduction}
Before we introduce Hua domains, we first recall the results on
generalized complex ellipsoids and bounded symmetric domains. A
generalized complex ellipsoid  (also called generalized
pseudoellipsoid) is a domain of the form
\begin{eqnarray*}
\begin{aligned}
\Sigma(\mathbf{n};\mathbf{p})&=\left\{(\zeta_1,\cdots,\zeta_r)\in
\mathbb{C}^{n_1}\times\cdots\times\mathbb{C}^{n_r}:
\sum_{k=1}^{r}\|\zeta_k\|^{2p_k}<1 \right\},
\end{aligned}
\end{eqnarray*}
where $\mathbf{n}=(n_1,\cdots,n_r)\in \mathbb{N}^{r},$
$\mathbf{p}=(p_1,\cdots,p_r)\in (\mathbb{R}_{+})^r,$ and $\|\cdot\|$
is the standard Hermitian norm. By relabelling the coordinates, we
can always  assume that $p_2\neq 1,\cdots,p_r\neq 1$, that is, there
is at most one 1 in $p_1,\cdots,p_r$.

In the special case where all the $p_k=1,$ the generalized complex
ellipsoid $\Sigma(\mathbf{n};\mathbf{p})$ reduces to the unit ball
in $\mathbb{C}^{n_1+\cdots+n_r}.$ Also, it is known that a
generalized complex ellipsoid $\Sigma(\mathbf{n};\mathbf{p})$ is
homogeneous if and only if $p_k=1$ for all $1\leq k\leq r$ (cf.
Kodama \cite{K1}). In general, a generalized complex ellipsoid is
not strongly pseudoconvex and its boundary is not smooth.

\vskip 6pt For the biholomorphic mappings between two
equidimensional  generalized complex ellipsoids, in 1968,  Naruki
\cite{Naruki} proved the following result.

\vskip 6pt \noindent {\bf Theorem 1.A} (Naruki \cite{Naruki}) {\it
Let $\Sigma(\mathbf{n};\mathbf{p})$ and
$\Sigma(\mathbf{m};\mathbf{q})$ be two equidimensional generalized
complex ellipsoids with $\mathbf{n},\; \mathbf{m}\in\mathbb{N}^r$
and $\mathbf{p},\; \mathbf{q}\in (\mathbb{R}_{+})^r$ (where $p_k\neq
1,\; q_k\neq 1$ for $2\leq k\leq r$). Then
$\Sigma(\mathbf{n};\mathbf{p})$ is biholomorphic to
$\Sigma(\mathbf{m};\mathbf{q})$ if and only if there exists a
permutation $\sigma\in S_r$ (where $S_r$ is the permutation group of
the $r$ numbers $\{1,\cdots,r\}$ ) such that $n_{\sigma(j)}=m_j,\;
p_{\sigma(j)}=q_j$ for $1\leq j\leq r$.}

\vskip 6pt

The holomorphic automorphism group ${\rm
Aut}(\Sigma(\mathbf{n};\mathbf{p}))$ of
$\Sigma(\mathbf{n};\mathbf{p})$ has been studied by Dini-Primicerio
\cite{Dini}, Kodama \cite{K1} and Kodama-Krantz-Ma \cite{Kodama}. In
2013, Kodama \cite{K1} obtained the result as follows.

\vskip 6pt \noindent {\bf Theorem 1.B} (Kodama \cite{K1}) {\it $(i)$
If $1$ does not appear in $p_1,\cdots,p_r$, then any automorphism
$\varphi\in {\rm Aut}(\Sigma(\mathbf{n};\mathbf{p}))$ is of the form
\begin{equation}
\varphi(\zeta_1,\cdots,\zeta_r)=\left(\gamma_1(\zeta_{\sigma(1)}),\cdots,\gamma_r(\zeta_{\sigma(r)})\right),
\end{equation}
where $\sigma\in S_r$ is a permutation of the $r$ numbers
$\{1,\cdots,r\}$ such that $n_{\sigma(i)}=n_i, p_{\sigma(i)}=p_i$
$(1\leq i\leq r)$ and $\gamma_1,\cdots, \gamma_r$ are unitary
transformations of $\mathbb{C}^{n_1}(n_{\sigma(1)}=n_1),\cdots,
\mathbb{C}^{n_r}(n_{\sigma(r)}=n_r)$ respectively.

$(ii)$ If $1$ appears in $p_1,\cdots,p_r$, we can assume, without
loss of generality, that $p_1=1, p_2\neq 1,\cdots,p_r\neq 1$, then
${\rm Aut}(\Sigma(\mathbf{n};\mathbf{p}))$ is generated by elements
of the form $(1)$ and automorphisms of the form
\begin{equation}
\varphi_a(\zeta_1,\zeta_2,\cdots,\zeta_r)=\left(T_{a}(\zeta_1),\zeta_2(\psi_a(\zeta_1))^{\frac{1}{2p_2}},
\cdots,\zeta_r(\psi_a(\zeta_1))^{\frac{1}{2p_r}}\right),
\end{equation}
where $T_a$ is an automorphism of the ball $\mathbf{B}^{n_1}$ in
$\mathbb{C}^{n_1}$, which brings a point $a\in \mathbf{B}^{n_1}$ in
the origin and
$$\psi_a(\zeta_1)=\frac{1-\|a\|^2}{(1-\langle
\zeta_1,a\rangle)^2}.$$}

\vskip 6pt Every bounded symmetric domain is, when equipped with the
Bergman metric, a Hermitian symmetric manifold of noncompact type,
and every Hermitian symmetric manifold of noncompact type can be
realized as a bounded symmetric domain in some $\mathbb{C}^d$ by the
Harish-Chandra embedding theorem. In 1935, E. Cartan proved that
there exist only six types of irreducible bounded symmetric domains.
They are four types of classical bounded symmetric domains and two
exceptional domains. So bounded symmetric domains are also known as
Cartan domains.

\vskip 6pt Let $\Omega$ be an irreducible bounded symmetric domain
in $\mathbb{C}^d$ of genus $g$ in its Harish-Chandra realization.
Let
$$\left\{\frac{1}{\sqrt{V(\Omega)}},h_1(z),h_2(z),\cdots\right\}$$
be an orthonormal basis of the Hilbert space $A^2(\Omega)$ of
square-integrable holomorphic functions on $\Omega$. Define the
Bergman kernel $K_{\Omega}(z,\bar{\xi})$ of $\Omega$ by
$$K_{\Omega}(z,\bar{\xi}):=\frac{1}{V(\Omega)}+\sum_{i=1}^{\infty}h_i(z)\overline{h_i(\xi)}$$
for all $z,\xi\in\Omega$. Obviously, $1\leq V(\Omega)K_{\Omega}(z,\bar{z})<+\infty$.
The generic norm of $\Omega$ is defined by
$$N_{\Omega}(z,\bar{\xi}):=\left(V(\Omega)K_{\Omega}(z,\bar{\xi})\right)^{-\frac{1}{g}}\;\;\; (z,\xi\in\Omega),$$
where $(V(\Omega)K_{\Omega}(z,\bar{\xi}))^{-\frac{1}{g}}:=\exp
(-\frac{1}{g}\log (V(\Omega)K_{\Omega}(z,\bar{\xi})))$, in which
$\log$ denotes the principal branch of logarithm (note
$K_{\Omega}(z,\bar{\xi})\neq 0$ for all $z,\xi\in \Omega$). Thus $0<
N_{\Omega}(z,\bar{z})\leq 1$ for all $z\in \Omega$ and
$N_{\Omega}(z,\bar{z})=0$ on the boundary of $\Omega$.

For an irreducible bounded symmetric domain $\Omega \subset
\mathbb{C}^d$ in its Harish-Chandra realization, a positive integer
$r$ and $\mathbf{n}=(n_1,\cdots,n_r)\in\mathbb{N}^r$,
$\mathbf{p}=(p_1,\cdots,p_r)\in(\mathbb{R}_{+})^r$, the Hua domain
${H_{\Omega}}(\mathbf{n};\mathbf{p})$ is defined by
\begin{align*}
&{H_{\Omega}}(\mathbf{n};\mathbf{p})={H_{\Omega}}(n_1,\cdots,n_r;p_1,\cdots,p_r)\\
&:=\left\{(z,w_{(1)},\cdots,w_{(r)})\in \Omega \times
\mathbb{C}^{n_1}\times\cdots\times\mathbb{C}^{n_r}:
\sum_{j=1}^r\|w_{(j)}\|^{2p_j}< N_{\Omega}(z,\bar{z})\right\},
\end{align*}
where $\|\cdot\|$ is the standard Hermitian norm. Note that
$\Omega\times \{0\}\subset {H_{\Omega}}(\mathbf{n};\mathbf{p})$ and
$b\Omega\times \{0\}\subset b{H}_{\Omega}(\mathbf{n};\mathbf{p})$
(where $bD$ denotes the boundary of a domain $D$).

For $(z,w_{(1)},\cdots,w_{(r)})\in
H_{\Omega}(\mathbf{n};\mathbf{p}),$ by definition, we have
\begin{align*}
\frac{\sum_{j=1}^r\|w_{(j)}\|^{2p_j}}{N_{\Omega}(z,\bar{z})}  =\exp
\left\{\log \left(\sum_{j=1}^r\|w_{(j)}\|^{2p_j}\right)+
{\frac{1}{g}} \log ( V(\Omega) K_{\Omega} (z,\bar{z}) )\right\}.
\end{align*}
Because $\log \left(\sum_{j=1}^r\|w_{(j)}\|^{2p_j}\right)+
{\frac{1}{g}} \log ( V(\Omega) K_{\Omega} (z,\bar{z}) )$ is a
 plurisubharmonic function on
$H_{\Omega}(\mathbf{n};\mathbf{p}),$ we have
${\sum_{j=1}^r\|w_{(j)}\|^{2p_j}}/{N_{\Omega}(z,\bar{z})}$ is a
continuous plurisubharmonic function on
$H_{\Omega}(\mathbf{n};\mathbf{p}).$  Since $\frac{1}{1-x}$ is a
monotonically increasing convex function for $x\in (-\infty,1)$ and
$0\leq \sum_{j=1}^r\|w_{(j)}\|^{2p_j}/N_{\Omega}(z,\bar{z})< 1$ on
$H_{\Omega}(\mathbf{n};\mathbf{p}),$ we have
$$\frac{1}{1-{\sum_{j=1}^r\|w_{(j)}\|^{2p_j}}/{N_{\Omega}(z,\bar{z})}}
=\frac{N_{\Omega}(z,\bar{z})}{N_{\Omega}(z,\bar{z})- \sum_{j=1}^r\|w_{(j)}\|^{2p_j}}$$ is a
continuous plurisubharmonic function on $H_{\Omega}(\mathbf{n};\mathbf{p}).$
Thus
$$\max\{\frac{N_{\Omega}(z,\bar{z})}{N_{\Omega}(z,\bar{z})- \sum_{j=1}^r\|w_{(j)}\|^{2p_j}}, \frac{1}{N_{\Omega}(z,\bar{z})}\}$$
is a continuous plurisubharmonic exhaustion function of $H_{\Omega}(\mathbf{n};\mathbf{p}).$
Then
$H_{\Omega}(\mathbf{n};\mathbf{p})$ is a bounded pseudoconvex domain
in $\mathbb{C}^{d+n_1+\cdots+n_r}.$ But, in general, a Hua domain is
a nonhomogeneous domain without smooth boundary.

Let $\mathcal{M}_{m,n}$ be the set of all $m\times n$ matrices
$z=(z_{ij})$ with complex entries. Let ${\overline z}$ be the
complex conjugate of the matrix $z$ and let ${z}^t$ be the transpose
of the matrix $z$. $I$ denotes the identity matrix. If a square
matrix $z$ is positive definite, then we write $z>0$. For each
bounded classical  symmetric domain $\Omega$ (refer to Hua
\cite{Hua}), we list the genus $g(\Omega)$, the generic norm
$N_\Omega(z,\overline{z})$ of $\Omega$ and corresponding Hua domain
${H_{\Omega}}(\mathbf{n};\mathbf{p})$ (see Yin-Wang-Zhao-Zhao-Guan
\cite{Yin1}) according to its type as following.

 $(i)$ If $\Omega=\Omega_I(m,n):=\{z\in
\mathcal{M}_{m,n}: I-z{\overline z}^t>0 \} \subset \mathbb{C}^d$
($1\leq m\leq n,\; d=mn$) (the classical domains of type $I$), then
$g(\Omega)=m+n,$ $N_\Omega(z,\overline{z})=\det(I-z{\overline
z}^t),$ and
$${H_{\Omega}}(\mathbf{n};\mathbf{p})=\left\{(z,w_{(1)},\cdots,w_{(r)})\in \Omega_I(m,n) \times
\mathbb{C}^{n_1}\times\cdots\times\mathbb{C}^{n_r}:
\sum_{j=1}^r\|w_{(j)}\|^{2p_j}< \det(I-z{\overline z}^t)\right\}.$$
Specially, when $\Omega=\mathbf{B}^d$ is the unit ball in $
\mathbb{C}^{d}$, then $N_\Omega(z,\overline{z})=1-\|z\|^2,$ and
$${H_{\Omega}}(\mathbf{n};\mathbf{p})=\left\{(z,w_{(1)},\cdots,w_{(r)})\in B^d \times
\mathbb{C}^{n_1}\times\cdots\times\mathbb{C}^{n_r}: \|z\|^2+
\sum_{j=1}^r\|w_{(j)}\|^{2p_j}< 1\right\}.$$ Thus the Hua domain
$H_{\mathbf{B}^d}(\mathbf{n};\mathbf{p})$ is just the generalized
complex ellipsoid $\Sigma((d,\mathbf{n});(1,\mathbf{p}))$.

$(ii)$ If $\Omega= \Omega_{II}(n):=\{z\in \mathcal{M}_{n,n}: z^t=-z,
I-z{\overline z}^t>0 \}\subset \mathbb{C}^d$ ($n\geq 2,\;
d={n(n-1)}/{2}$) (the classical domains of type $II$), then
$g(\Omega)=2(n-1),$ $N_\Omega(z,\overline{z})=(\det(I-z{\overline
z}^t))^{1/2} ,$ and
$${H_{\Omega}}(\mathbf{n};\mathbf{p})=\left\{(z,w_{(1)},\cdots,w_{(r)})\in \Omega_{II}(n) \times
\mathbb{C}^{n_1}\times\cdots\times\mathbb{C}^{n_r}:
\sum_{j=1}^r\|w_{(j)}\|^{2p_j}< (\det(I-z{\overline
z}^t))^{1/2}\right\}.$$

$(iii)$ If $\Omega= \Omega_{III}(n):=\{z\in \mathcal{M}_{n,n}:
z^t=z, I-z{\overline z}^t>0 \}\subset \mathbb{C}^d$ ($n\geq 2,\;
d=n(n+1)/2$) (the classical domains of type $III$), then
$g(\Omega)=n+1,$ $N_\Omega(z,\overline{z})=\det(I-z{\overline z}^t)
,$ and
$${H_{\Omega}}(\mathbf{n};\mathbf{p})=\left\{(z,w_{(1)},\cdots,w_{(r)})\in \Omega_{III}(n) \times
\mathbb{C}^{n_1}\times\cdots\times\mathbb{C}^{n_r}:
\sum_{j=1}^r\|w_{(j)}\|^{2p_j}< \det(I-z{\overline z}^t)\right\}.$$

$(iv)$ If $\Omega=\Omega_{IV}(n):=\{z\in \mathbb{C}^{n}:
1-2z{\overline z}^t+|zz^t|^2>0,  z{\overline z}^t<1\} $ ($n\geq 3$)
(the classical domains of type $IV$), then
 $g(\Omega)=n,$
$N_\Omega(z,\overline{z})=1-2z{\overline z}^t +|zz^t|^2,$ and
$${H_{\Omega}}(\mathbf{n};\mathbf{p})=\left\{(z,w_{(1)},\cdots,w_{(r)})\in \Omega_{IV}(n) \times
\mathbb{C}^{n_1}\times\cdots\times\mathbb{C}^{n_r}:
\sum_{j=1}^r\|w_{(j)}\|^{2p_j}<1-2z{\overline z}^t
+|zz^t|^2\right\}.$$

Let $\Omega$ be an irreducible bounded symmetric domain in
$\mathbb{C}^d$ in its Harish-Chandra realization. We can always
assume that the Hua domain ${H_{\Omega}}(\mathbf{n},\mathbf{p})$ is
written {\it in its standard form}, that is,

(i) If $\Omega$ is the unit ball, then $p_1\neq 1,\cdots, p_r\neq 1$
 (here it is understood that this domain is the unit ball in
$\mathbb{C}^d$ if $r=0.$);

(ii) If ${\rm rank}(\Omega)\geq 2$, then $p_1=1,p_2\neq
1,\cdots,p_r\neq 1$ (here it is understood that $p_1=1$ does not
appear if $n_1=0$).

It is easy to see that every Hua domain can be written {\it in its
standard form} by relabelling the coordinates. Therefore, for every
given Hua domain, there exists an irreducible bounded symmetric
domain $\Omega$ in its Harish-Chandra realization such that the Hua
domain can be written as ${H_{\Omega}}(\mathbf{n},\mathbf{p})$ {\it
in its standard form}.

Let $\Omega$ be an irreducible bounded symmetric domain in
$\mathbb{C}^d$, and $\mathbf{n}\in\mathbb{N}^r$,
$\mathbf{p}\in(\mathbb{R}_{+})^r$. Let the family
$\Gamma({H_{\Omega}}(\mathbf{n};\mathbf{p}))$ be exactly the set of
all mappings $\Phi$:
\begin{equation}
\Phi(z,w_{(1)},\cdots,w_{(r)}) =\left(\varphi(z),\;
U_1(w_{(1)})\frac{
(N_{\Omega}(z_0,\overline{z_0}))^{\frac{1}{2p_1}}}{(N_{\Omega}(z,\overline{z_0}))^{\frac{1}{p_1}}},\;
\cdots,\;
U_r(w_{(r)})\frac{(N_{\Omega}(z_0,\overline{z_0}))^{\frac{1}{2p_r}}}{(N_{\Omega}(z,\overline{z_0}))^{\frac{1}{p_r}}}\right)
\end{equation}
for $(z,w_{(1)},\cdots,w_{(r)})\in
{H_{\Omega}}(\mathbf{n};\mathbf{p})$, where $\varphi\in {\rm
Aut}(\Omega)$, $U_j$ is a unitary transformation of
$\mathbb{C}^{n_j}$ for $1\leq j\leq r$, and $z_0=\varphi^{-1}(0)$.
Then $\Gamma ({H_{\Omega}}(\mathbf{n};\mathbf{p}))$ is a subgroup of
the holomorphic automorphism group ${\rm
Aut}({H_{\Omega}}(\mathbf{n};\mathbf{p}))$ of
${H_{\Omega}}(\mathbf{n};\mathbf{p})$ (see Yin-Wang-Zhao-Zhao-Guan
\cite{Yin1}). Obviously, every element of $\Gamma
({H_{\Omega}}(\mathbf{n};\mathbf{p}))$ preserves the set
$\Omega\times \{0\}(\subset {H_{\Omega}}(\mathbf{n};\mathbf{p}))$
and $\Gamma ({H_{\Omega}}(\mathbf{n};\mathbf{p}))$ is transitive on
$\Omega\times \{0\}(\subset {H_{\Omega}}(\mathbf{n};\mathbf{p}))$.
For the general reference of Hua domains, see
Yin-Wang-Zhao-Zhao-Guan \cite{Yin1} and references therein.

When $r=1$, the  Hua domain ${H_{\Omega}}(n_1;p_1)$ is also called
the Cartan-Hartogs domain and is also denoted by
$\Omega^{B^{n_1}}(p_1)$. For the reference of the Cartan-Hartogs
domains, see Ahn-Byun-Park \cite{ABP},  Feng-Tu \cite{Tu3},
Loi-Zedda \cite{LZ}, Wang-Yin-Zhang-Roos \cite{RWYZ} and Yin
\cite{Yin} and references therein.

In 2012, Ahn-Byun-Park \cite{ABP} determined the automorphism group
of the Cartan-Hartogs domain ${H_{\Omega}}(n_1;p_1)$ by case-by-case
checking only for four types of classical domains $\Omega$.
Following the reasoning in Ahn-Byun-Park \cite{ABP}, Rong \cite{R}
claimed a description of automorphism groups of Hua domains
${H_{\Omega}}(\mathbf{n},\mathbf{p})$ in 2014. But, Lemma 3.2 in
Rong \cite{R}, which is central to the proof of its main results in
\cite{R}, is definitely wrong (cf. Proposition 2.4 in our paper for
references).

\vskip 6pt The first goal of this paper is to give a description of
the biholomorphisms between two equidimensional Hua domains. By
using a different technique  from that in  Ahn-Byun-Park \cite{ABP},
we obtain the result as follows.

\vskip 6pt \noindent {\bf Theorem 1.1}. {\it Suppose that
$$f:H_{\Omega_1}(\mathbf{n};\mathbf{p})\rightarrow H_{\Omega_2}(\mathbf{m};\mathbf{q})$$
is a biholomorphism between two equidimensional Hua domains
$H_{\Omega_1}(\mathbf{n};\mathbf{p})$ and
$H_{\Omega_2}(\mathbf{m};\mathbf{q})$ in their standard forms, where
$\Omega_1\subset \mathbb{C}^{d_1}$ and $\Omega_2\subset
\mathbb{C}^{d_2}$ are two irreducible bounded symmetric domains in
the Harish-Chandra realization, and $\mathbf{n},
\mathbf{m}\in\mathbb{N}^r$,
$\mathbf{p},\mathbf{q}\in(\mathbb{R}_{+})^r$. Then there exists an
automorphism $\Phi\in \Gamma (H_{\Omega_2}(\mathbf{m};\mathbf{q}))$
(see $(3)$ here) and a permutation $\sigma\in S_r$ with
$n_{\sigma(i)}=m_i, p_{\sigma(i)}=q_i$ for $1\leq i\leq r$ such that
\begin{eqnarray*}
\begin{aligned}
\Phi\circ f(z,w_{(1)},\cdots,w_{(r)})=(z,w_{(\sigma(1))},\cdots,w_{(\sigma(r))})
\begin{pmatrix}A &   &  &  \\
  & U_1 &  &  \\
  &     &\ddots &  \\
  &     &  & U_r\\
\end{pmatrix},
\end{aligned}
\end{eqnarray*}
where $A$ is a complex linear isomorphism of
$\mathbb{C}^d\;(d:=d_1=d_2)$ with $A(\Omega_1)=\Omega_2$, and $U_i$
is a unitary transformation of $\mathbb{C}^{m_i}\;
(m_i=n_{\sigma(i)})$ for $1\leq i\leq r$.}

\vskip 6pt As a special result of Theorem 1.1, we completely
describe the automorphism group of the Hua domains
${H_{\Omega}}(\mathbf{n};\mathbf{p})$ for all irreducible bounded
symmetric domains $\Omega$ as follows.

\vskip 6pt \noindent {\bf Corollary 1.2}. {\it Let
${H_{\Omega}}(\mathbf{n};\mathbf{p})$ be a Hua domain in its
standard form and $\Gamma({H_{\Omega}}(\mathbf{n};\mathbf{p}))$ is
generated by the mappings of the form $(3)$, where $\Omega\subset
\mathbb{C}^{d}$ is an irreducible bounded symmetric domain in the
Harish-Chandra realization, and $\mathbf{n}\in\mathbb{N}^r$,
$\mathbf{p}\in(\mathbb{R}_{+})^r$. Then, for every $f\in {\rm
Aut}({H_{\Omega}}(\mathbf{n};\mathbf{p})),$ there exist a $\Phi\in
\Gamma({H_{\Omega}}(\mathbf{n};\mathbf{p}))$ and a permutation
$\sigma\in S_r$ with $n_{\sigma(i)}=n_i,\; p_{\sigma(i)}=p_i$ for
$1\leq i\leq r$ such that
$$f(z,w_{(1)},\cdots,w_{(r)})=\Phi(z,w_{(\sigma(1))},\cdots,w_{(\sigma(r))}).$$}

Remarks on Corollary 1.2. In the case of $r=1$, Ahn-Byun-Park
\cite{ABP} obtained Corollary 1.2 for four types of classical
domains $\Omega$ in 2012. If $\Omega=\mathbf{B}^d$, then the Hua
domain
$H_{\mathbf{B}^d}(\mathbf{n},\mathbf{p})=\Sigma((d,\mathbf{n}),(1,\mathbf{p}))$
is a generalized complex ellipsoid. By (3), we have that Corollary
1.2 implies Theorem 1.B (ii).

It is important that the Hua domain is written {\it in its standard
form} in Corollary 1.2. (i) For example, define
$$H_{\mathbf{B}^2}(\mathrm{(2,2)};\mathrm{(1,2)})=\left\{(z,w_{(1)},w_{(2)} )\in \mathbf{B}^2 \times
\mathbb{C}^2 \times  \mathbb{C}^2: \|z\|^2+ \|w_{(1)}\|^{2} +
\|w_{(2)}\|^{4}< 1\right\}.$$ Then, in this case, there exists only
the identity $\sigma=1\in S_2$ such that $n_{\sigma(i)}=n_i,\;
p_{\sigma(i)}=p_i$ for $1\leq i\leq 2$, and obviously
$\Gamma(H_{\mathbf{B}^2}(\mathrm{(2,2)};\mathrm{(1,2)}))\subsetneqq
{\rm Aut}(H_{\mathbf{B}^2}(\mathrm{(2,2)};\mathrm{(1,2)}))$ (cf.
Theorem 1.1 in Rong \cite{R}). This means that Corollary 1.2 does
not hold for the Hua domain ${H_{\Omega}}(\mathbf{n};\mathbf{p})$
which is not in the standard form. (ii) If the Hua domain
${H_{\Omega}}(\mathbf{n};\mathbf{p})$ is the unit ball, then we have
${H_{\Omega}}(\mathbf{n};\mathbf{p})=\Omega$ and $r=0$ (by
${H_{\Omega}}(\mathbf{n};\mathbf{p})$ in its standard form).
Therefore, we have
$$\Gamma({H_{\Omega}}(\mathbf{n};\mathbf{p})) ={\rm Aut}(\Omega)\; (=
{\rm Aut}({H_{\Omega}}(\mathbf{n};\mathbf{p}))).$$ This means that
Corollary 1.2 holds for the unit ball case of Hua domain
${H_{\Omega}}(\mathbf{n};\mathbf{p})$  in its standard form (cf.
Theorem 1.1 in Ahn-Byun-Park \cite{ABP}).

\vskip 10pt The second purpose of this paper is to study proper
holomorphic mappings between Hua domains. We first recall the
structure of proper holomorphic self-mappings of the unit ball
$\mathbf{B}^n$ in $\mathbb{C}^n$. When $n=1$, such maps are
precisely the finite Blaschke products. The situation is quite
different for $n\geq 2$. The following fundamental result was proved
by Alexander \cite{Ale} in 1977.

\vskip 6pt \noindent {\bf Theorem 1.C} (Alexander \cite{Ale}) {\it
Any proper holomorphic self-mapping of the unit ball $\mathbf{B}^n$
in $\mathbb{C}^{n}$ $(n\geq 2)$ is an automorphism of
$\mathbf{B}^n$.}

\vskip 6pt We remark that $$f(z_1,z_2)=(z_1,z_2^2):\;
|z_1|^2+|z_2|^4<1 \longrightarrow |w_1|^2+|w_2|^2<1$$ is a proper
holomorphic mapping between two bounded pseudoconvex domains in
$\mathbb{C}^2$ with smooth real-analytic boundary, but it is
branched and is not biholomorphic. Thus it suggests a subject to
discover some interesting bounded weakly pseudoconvex domains
$D_1,\;D_2$ in {\bf{C}}$^n$ $(n\geq 2)$ such that any proper
holomorphic mapping from $D_1$ to $D_2$ is a biholomorphism. There
are many important results concerning proper holomorphic mapping $f:
D_1\rightarrow D_2$ between two bounded pseudoconvex domains $D_1,\;
D_2$ in $\mathbb{C}^n$ with smooth boundary. If the proper
holomorphic mapping $f$ extends smoothly to the closure of $D_1$,
then the extended mapping takes the boundary $bD_1$ into the
boundary $bD_2$, and it satisfies the tangential Cauchy-Riemann
equations on $bD_1$. Thus the proper holomorphic mapping
$f:D_1\rightarrow D_2$ leads naturally to the geometric study of the
mappings from $bD_1$ into $bD_2$. These researches are often heavily
based on analytic techniques about the mapping on boundaries (e.g.,
see Forstneri\v{c} \cite{Forst}  and Huang \cite{H}). The lack of
boundary regularity usually presents a serious analytical
difficulty.

As we know, in general, a generalized complex ellipsoid is not
strongly pseudoconvex and its boundary is not smooth. Also, there
are many results (e.g., Dini-Primicerio \cite{Dini91, Dini}, Hamada
\cite{H} and  Landucci \cite{L}) concerning proper holomorphic
mappings between two generalized complex ellipsoids.

\vskip 6pt For the case of $\mathbf{p},\mathbf{q}\in
(\mathbb{Z}_+)^r$, in 1997, Dini-Primicerio (\cite{Dini}, Th. 4.6)
proved the following result.

\vskip 6pt \noindent {\bf Theorem 1.D} (Dini-Primicerio \cite{Dini})
{\it Let $\Sigma(\mathbf{n};\mathbf{p})$ and
$\Sigma(\mathbf{m};\mathbf{q})$ be two equidimensional generalized
complex ellipsoids with $\mathbf{n},\mathbf{m}\in\mathbb{N}^r$ and
$\mathbf{p},\mathbf{q}\in (\mathbb{Z}_+)^r$ (where $p_k\neq 1,\;
q_k\neq 1$ for $2\leq k\leq r$) such that $n_k\geq 2$ whenever
$p_k\geq 2$  and $m_k\geq 2$ whenever $q_k\geq 2$  for $1\leq k\leq
r$. Then there exists a proper holomorphic mapping
$f:\Sigma(\mathbf{n};\mathbf{p})\rightarrow
\Sigma(\mathbf{m};\mathbf{q})$ if and only if there exists a
permutation $\sigma\in S_r$ such that $n_{\sigma(j)}=m_j,\;
p_{\sigma(j)}=q_j$ for $1\leq j\leq r$.}

\vskip 6pt

Remark. When $\mathbf{p},\mathbf{q}\in (\mathbb{Z}_+)^r$, we have
that $\Sigma(\mathbf{n};\mathbf{p})$ and
$\Sigma(\mathbf{m};\mathbf{q})$ are pseudoconvex domains with real
analytic boundaries. Theorem 1.D comes from Theorem 4.6 in
Dini-Primicerio \cite{Dini}. In Dini-Primicerio \cite{Dini}, Theorem
4.6 is proved by Theorem 3.1  (in Dini-Primicerio \cite{Dini})
assuming that ``the sets of weak pseudoconvexity of
$\Sigma(\mathbf{n};\mathbf{p})$ and $\Sigma(\mathbf{m};\mathbf{q})$
are contained in analytic sets of codimension at least 2", which is
equivalent to ``$n_k\geq 2$ whenever $p_k\geq 2$  and $m_k\geq 2$
whenever  $q_k\geq 2$  for $1\leq k\leq r$" in Theorem 1.D (see
(2.2) in \cite{Dini} for references).

\vskip 6pt Following the methods of Pinchuk \cite{Pin},
Dini-Primicerio \cite{Dini} proved the so called ``localization
principle of biholomorphisms" for generalized complex ellipsoids,
that is, any local biholomorphism sending boundary points to
boundary points extends to a global one, and, as its application,
Dini-Primicerio \cite{Dini} get Theorem 1.D. The approach of
Pinchuck \cite{Pin} to ``localization principle of biholomorphisms"
is firstly to show that the local biholomorphism is rational (thus
extends naturally to be globally meromorphic),  and then to show
that the rational mapping is biholomorphic by the standard argument:
if the zero locus of the holomorphic Jacobian  determinant of the
rational mapping is nonempty, then the set of points of weak
pseudoconvexity should contain a set of real codimension $3$. Thus
the assumption that ``the sets of weak pseudoconvexity is contained
in some complex analytic set of complex codimension at least 2" will
force the zero locus of the  holomorphic Jacobian  determinant to be
empty. Thus the conditions
``$\mathbf{p},\mathbf{q}\in(\mathbb{Z}_+)^r$" and ``$n_i\geq 2$
whenever $p_i\geq 2$  and $m_i\geq 2$  whenever  $q_i\geq 2$  for
$1\leq i\leq r$" are indispensable in proving Theorem 1.D.

\vskip 6pt Even though the bounded homogeneous domains in
$\mathbb{C}^n$ are always pseudoconvex, there are, of course, many
such domains (e.g., all bounded symmetric domains of rank $\geq 2$)
such that they do not have smooth boundary and have no strongly
pseudoconvex boundary point by the Wong-Rosay theorem (see Rudin
\cite{R1}, Theorem 15.5.10 and its Corollary). There are many
rigidity results about the proper holomorphic mappings between
bounded symmetric domains.

In 1984, by using results of Bell \cite{Bell} and Tumanov-Henkin
\cite{Tuma}, Henkin-Novikov \cite{Henkin} proved the following
result (see Th.3.3 in Forstneri\v{c} \cite{Forst} for references).

\vskip 6pt \noindent {\bf Theorem 1.E} (Henkin-Novikov
\cite{Henkin}) {\it Any proper holomorphic self-mapping on an
irreducible bounded symmetric domain of rank $\geq 2$ is an analytic
automorphism.}

\vskip 6pt

Using the idea in Mok-Tsai \cite{Mok-Tsai} and Tsai \cite{Tsai}, Tu
\cite{Tu1,Tu2} (one of the authors of the current article)
obtained rigidity results on proper holomorphic mappings between
bounded symmetric domains and proved the following in 2002.

\vskip 6pt \noindent {\bf Theorem 1.F} (Tu \cite{Tu1}) {\it Let
$\Omega_1$ and $\Omega_2$ be two equidimensional bounded symmetric
domains. Assume that $\Omega_1$ is irreducible and
$rank(\Omega_1)\geq 2$. Then, any proper holomorphic mapping from
$\Omega_1$ to $\Omega_2$ is a biholomorphism.}

\vskip 6pt

Further, using the idea in Mok-Tsai \cite{Mok-Tsai} and Tsai
\cite{Tsai}, in 2010, Mok-Ng-Tu \cite{MNT} obtained some rigidity
results of proper holomorphic mappings on bounded symmetric domains
as follows.

\vskip 6pt \noindent
 {\bf  Theorem 1.G} (Mok-Ng-Tu \cite{MNT}) {\it Let  $\Omega$  be an
irreducible bounded symmetric domain of rank $\geq 2$ which is not
isomorphic to a Type-IV classical symmetric domain $D_{IV}^N$ of
dimension  $N \geq 3$. Let $F : \Omega\rightarrow  D$ be a proper
holomorphic map onto a bounded convex domain $D$. Then, $F :
\Omega\rightarrow  D$ is a biholomorphism and $D$ is, up to an
affine-linear transformation, the Harish-Chandra realization of
$\Omega$.}

\vskip 6pt The second goal of this paper is to establish what seems
to be the first rigidity result for proper holomorphic mappings on
Hua domains.

For a Hua domain
${H_{\Omega}}(\mathbf{n};\mathbf{p})={H_{\Omega}}(n_1,\cdots,n_r;p_1,\cdots,p_r)$
in its standard form. The boundary
$bH_{\Omega}(\mathbf{n};\mathbf{p})$ of
$H_{\Omega}(\mathbf{n};\mathbf{p})$  is comprised of
\begin{equation}\label{1.4}
bH_{\Omega}(\mathbf{n};\mathbf{p})
=b_0{H_{\Omega}}(\mathbf{n};\mathbf{p})\cup
b_1{H_{\Omega}}(\mathbf{n};\mathbf{p})\cup (b\Omega\times \{0\}),
\end{equation}
where
\begin{align*}
b_0{H_{\Omega}}(\mathbf{n};\mathbf{p}):=\Big \{&
(z,w_{(1)},\cdots,w_{(r)})\in\Omega\times
\mathbb{C}^{n_1}\times\cdots\times\mathbb{C}^{n_r}:\\
& \sum_{i=1}^r\|w_{(i)}\|^{2p_i}=N_{\Omega}(z,z),\|w_{(j)}\|^2\neq
0, 1+\delta\leq j\leq r\Big\},
\end{align*}
\begin{align*}
b_1{H_{\Omega}}(\mathbf{n};\mathbf{p}):= \bigcup_{j=1+\delta}^r
\Big\{& (z,w_{(1)},\cdots,w_{(r)})
\in\Omega\times \mathbb{C}^{n_1}\times\cdots\times\mathbb{C}^{n_r}:\\
&\sum_{i=1}^r\|w_{(i)}\|^{2p_i}=N_{\Omega}(z,z),\|w_{(j)}\|^2=0\Big\},
\end{align*}
in which
\begin{equation*}
\delta=
\begin{cases}
1 & \text{if $p_1=1$}, \\
0 & \text{if $p_1\neq 1$}.
\end{cases}
\end{equation*}
Then we have (by Proposition 2.4 in this paper):

$(a)$ $b_0{H_{\Omega}}(\mathbf{n};\mathbf{p})$ is a real analytic
hypersurface in $\mathbb{C}^{d+|\mathbf{n}|}$ and
${H_{\Omega}}(\mathbf{n};\mathbf{p})$ is strongly pseudoconvex at
all points of $b_0{H_{\Omega}}(\mathbf{n};\mathbf{p})$.

$(b)$ If ${H_{\Omega}}(\mathbf{n};\mathbf{p})$ isn't a ball, then
${H_{\Omega}}(\mathbf{n};\mathbf{p})$ is not strongly pseudoconvex
at any point of $b_1{H_{\Omega}}(\mathbf{n};\mathbf{p})\cup
(b\Omega\times \{0\}).$

Obviously, $b_1{H_{\Omega}}(\mathbf{n};\mathbf{p})\cup
(b\Omega\times \{0\})$ is contained in a complex analytic subset in
$\mathbb{C}^{d+|\mathbf{n}|}$ of complex codimension
$\min\{n_{1+\delta},\cdots,n_r, n_1+\cdots +n_r\}$ (note
$\min\{n_{1+\delta},\cdots,n_r, n_1+\cdots +n_r\}=n_1$ for $r=1$ and
$\min\{n_{1+\delta},\cdots,n_r, n_1+\cdots +n_r\}
=\min\{n_{1+\delta},\cdots,n_r\}$ for $r\geq 2$).

\vskip 6pt \noindent {\bf Theorem 1.3}. {\it Suppose that
$$f:H_{\Omega_1}(\mathbf{n}_1;\mathbf{p}_1)\rightarrow H_{\Omega_2}(\mathbf{n}_2;\mathbf{p}_2)$$
is a proper holomorphic mapping between two equidimensional Hua
domains $H_{\Omega_1}(\mathbf{n}_1;\mathbf{p}_1)$ and
$H_{\Omega_2}(\mathbf{n}_2;\mathbf{p}_2)$ in their standard forms,
where $\Omega_1\subset \mathbb{C}^{d_1}$ and $\Omega_2\subset
\mathbb{C}^{d_2}$ are two irreducible bounded symmetric domains in
the Harish-Chandra realization, and
$\mathbf{n}_1,\mathbf{n}_2\in\mathbb{N}^r$,
$\mathbf{p}_1,\mathbf{p}_2\in(\mathbb{R}_{+})^r$. Assume that
$b_1H_{\Omega_i}(\mathbf{n}_i;\mathbf{p}_i) \cup (b\Omega_i\times
\{0\})$  $(i=1,2)$ is contained in some complex analytic set of
complex codimension at least $2$. Then
$f:H_{\Omega_1}(\mathbf{n}_1;\mathbf{p}_1)\rightarrow
H_{\Omega_2}(\mathbf{n}_2;\mathbf{p}_2)$ is a biholomorphism. }

\vskip 6pt Remarks on Theorem 1.3. (i) In Theorem 1.3, we don't
assume $\dim \Omega_1=\dim\Omega_2$.

(ii) In Theorem 1.3,  the assumption
``$b_1H_{\Omega}(\mathbf{n};\mathbf{p}) \cup (b\Omega\times \{0\})$
is contained in some complex analytic set of complex codimension at
least $2$" is equivalent to that $H_{\Omega}(\mathbf{n};\mathbf{p})$
(in its standard form) satisfies
$$\min\{n_{1+\delta},\cdots,n_r, n_1+\cdots +n_r\}\geq 2,$$
that is, $H_{\Omega}(\mathbf{n};\mathbf{p})$ (in its standard form)
satisfies the following assumptions:
 (a) If $\Omega=\mathbf{B}_{d}$ is the unit ball, then
$\min\{n_{1},\cdots,n_r\}\geq 2;$ (b) If ${\rm rank}(\Omega)\geq 2$
and $p_1\neq 1$, then $\min\{n_{1},\cdots,n_r\}\geq 2;$ (c) If ${\rm
rank}(\Omega)\geq 2$ and $p_1=1$, then $\min\{n_{2},\cdots,n_r,n_1+n_2+\cdots+n_r\}\geq
2.$

(iii) In Theorem 1.3, the assumption ``
$b_1H_{\Omega_i}(\mathbf{n}_i;\mathbf{p}_i) \cup (b\Omega_i\times
\{0\})$  $(i=1,2)$ is contained in some complex analytic set of
complex codimension at least $2$" cannot be removed. For example,
let $\Omega$ be an irreducible bounded symmetric domain with ${\rm
rank}(\Omega)\geq 2$, $n_1:=1$ (i.e., $w_{(1)}\in \mathbb{C}$ ), and
$$\Phi(z,w_{(1)},w_{(2)},\cdots,w_{(r)}):=(z,w_{(1)}^2,w_{(2)},\cdots,w_{(r)})$$
for $(z,w_{(1)},w_{(2)},\cdots,w_{(r)})\in
{H_{\Omega}}(1,n_2,\cdots,n_r;p_1,p_2,\cdots,p_r)$. Then $\Phi$ is a
proper holomorphic mapping from ${H_{\Omega}}(1,n_2,\cdots,
n_r;{p_1},p_2,\cdots,p_r)$ to ${H_{\Omega}}(1,n_2,\cdots,n_r;
p_1/2,p_2,\cdots,p_r)$, but $\Phi$ is not a biholomorphism.

\vskip 6pt Combining Theorem 1.3 and Corollary 1.2, we immediately
have the result as follows.

\vskip 6pt

\noindent {\bf Corollary 1.4}. {\it Suppose that $f$
 is a proper holomorphic self-mapping on the Hua domain
 ${H_{\Omega}}(\mathbf{n};\mathbf{p})$
in its standard form, where $\Omega\subset \mathbb{C}^{d}$ is an
irreducible bounded symmetric domain in the Harish-Chandra
realization, and $\mathbf{n}\in\mathbb{N}^r$,
$\mathbf{p}\in(\mathbb{R}_{+})^r$  with
$\min\{n_{1+\delta},\cdots,n_r,n_1+n_2+\cdots+n_r\}\geq 2$. Then $f$
is an automorhism of the Hua domain
 ${H_{\Omega}}(\mathbf{n};\mathbf{p})$, that is,
there exist a $\Phi\in \Gamma({H_{\Omega}}(\mathbf{n};\mathbf{p}))$
and a permutation $\sigma\in S_r$ with $n_{\sigma(i)}=n_i,\;
p_{\sigma(i)}=p_i$ for $1\leq i\leq r$ such that
$$f(z,w_{(1)},\cdots,w_{(r)})=\Phi(z,w_{(\sigma(1))},\cdots,w_{(\sigma(r))}).$$}

When $\Omega\subset \mathbb{C}^d$ is the unit ball $\mathbf{B}^d$,
we get that
$H_{\mathbf{B}^d}(\mathbf{n};\mathbf{p})=\Sigma((d,\mathbf{n});(1,\mathbf{p}))$
is a generalized complex ellipsoid. Thus, by Theorem 1.3 and Theorem
1.1, we get the following result about proper holomorphic mappings
between generalized complex ellipsoids.

\vskip 6pt \noindent {\bf Corollary 1.5}. {\it Let
$\Sigma(\mathbf{n};\mathbf{p})$ and $\Sigma(\mathbf{m};\mathbf{q})$
be two equidimensional generalized complex ellipsoids with
$\mathbf{n},\mathbf{m}\in \mathbb{N}^r$ and
$\mathbf{p},\mathbf{q}\in (\mathbb{R}_+)^r$  (where $p_k\neq 1,\;
q_k\neq 1$ for $2\leq k\leq r$). Assume that $n_i\geq 2,$  $m_i\geq
2$  for $2\leq i\leq r$ and $p_1=1,$ $q_1=1$. Then there exists a
proper holomorphic mapping
$f:\Sigma(\mathbf{n};\mathbf{p})\rightarrow
\Sigma(\mathbf{m};\mathbf{q})$ if and only if there exists a
permutation $\sigma\in S_r$ such that $n_{\sigma(j)}=m_j,\;
p_{\sigma(j)}=q_j$ for $1\leq j\leq r$.}

\vskip 6pt Remark. When $\mathbf{p},\mathbf{q}\in (\mathbb{R}_+)^r$,
we have that, in general, $\Sigma(\mathbf{n};\mathbf{p})$ and
$\Sigma(\mathbf{m};\mathbf{q})$ are pseudoconvex domains without
smooth boundaries. Corollary 1.5 is an extension of Theorem 1.D to
the special case of $\mathbf{p},\mathbf{q}\in(\mathbb{R}_{+})^r$.

\vskip 6pt Now we shall present an outline of the argument in our
proof of main results.

In general, a Hua domain is a nonhomogeneous domain without smooth
boundary. But it is still a bounded complete circular domain. Let
$$f:\; H_{\Omega_1}(\mathbf{n};\mathbf{p})\rightarrow
H_{\Omega_2}(\mathbf{m},\mathbf{q})$$ be a proper holomorphic
mapping between two equidimensional Hua domains in their standard
forms. We want to prove that $f$ is a biholomorphism, and further,
to determine the explicit form of the biholomorphism $f$.

In order to prove that $f:\;
H_{\Omega_1}(\mathbf{n};\mathbf{p})\rightarrow
H_{\Omega_2}(\mathbf{m},\mathbf{q})$ is a biholomorphism, it
suffices to show that $f:\;
H_{\Omega_1}(\mathbf{n};\mathbf{p})\rightarrow
H_{\Omega_2}(\mathbf{m},\mathbf{q})$ is unbranched. The
transformation rule for Bergman kernels under proper holomorphic
mapping (e.g., Th. 1 in Bell \cite{Bell82}) plays a key role in
extending proper holomorphic mapping. Our idea here is heavily based
on the framework of Bell \cite{Bell, Bell82} and Pin\v{c}uk
\cite{Pin}. The first is to prove that $f$ extends holomorphically
to the closures. By using a kind of semi-regularity at the boundary
of the Bergman kernel associated to a Hua domain, we get the
extension by using the standard argument in Bell \cite{Bell}. The
second is to prove that $f:\;
H_{\Omega_1}(\mathbf{n};\mathbf{p})\rightarrow
H_{\Omega_2}(\mathbf{m},\mathbf{q})$ is unbranched assuming the
first one is achieved. By investigating the strongly pseudoconvex
part of the boundary of the Hua domains and using the local
regularity for the mappings between strongly pseudoconvex
hypersurfaces (e.g., see Pin\v{c}uk \cite{Pin}), we get that $f:\;
H_{\Omega_1}(\mathbf{n};\mathbf{p})\rightarrow
H_{\Omega_2}(\mathbf{m},\mathbf{q})$ is unbranched. So $f:\;
H_{\Omega_1}(\mathbf{n};\mathbf{p})\rightarrow
H_{\Omega_2}(\mathbf{m},\mathbf{q})$ is a biholomorphism.
Furthermore, by the uniqueness theorem, we have $$f:\;
H_{\Omega_1}(\mathbf{n};\mathbf{p})\rightarrow
H_{\Omega_2}(\mathbf{m},\mathbf{q})$$ extends to a biholomorphism
between their closures.

Next we show that $f:\;
H_{\Omega_1}(\mathbf{n};\mathbf{p})\rightarrow
H_{\Omega_2}(\mathbf{m},\mathbf{q})$ maps the base space to the base
space (that is, $f(\Omega_1\times \{0\})\subset
\Omega_2\times\{0\}$). Let
$$bH_{\Omega_1}(\mathbf{n};\mathbf{p}):=
b_0{H_{\Omega_1}}(\mathbf{n};\mathbf{p})\cup
b_1{H_{\Omega_1}}(\mathbf{n};\mathbf{p})\cup (b\Omega_1\times
\{0\}),$$
$$ bH_{\Omega_2}(\mathbf{m};\mathbf{q}):=
b_0{H_{\Omega_2}}(\mathbf{m};\mathbf{q})\cup
b_1{H_{\Omega_2}}(\mathbf{m};\mathbf{q})\cup (b\Omega_2\times
\{0\}),
$$
where see $(4)$ for the notations. Then $(a)$
$H_{\Omega_1}(\mathbf{n};\mathbf{p})$ (resp.,
$H_{\Omega_2}(\mathbf{m};\mathbf{q})$) is strongly pseudoconvex at
all points of $b_0{H_{\Omega_1}}(\mathbf{n};\mathbf{p})$ (resp.,
$b_0H_{\Omega_2}(\mathbf{m};\mathbf{q})$); $(b)$ If
${H_{\Omega_1}}(\mathbf{n};\mathbf{p})$ (resp.,
$H_{\Omega_2}(\mathbf{m};\mathbf{q})$) isn't a ball, then
$H_{\Omega_1}(\mathbf{n};\mathbf{p})$ (resp.,
$H_{\Omega_2}(\mathbf{m};\mathbf{q})$) is not strongly pseudoconvex
at any point of $b_1{H_{\Omega_1}}(\mathbf{n};\mathbf{p})\cup
(b\Omega_1\times \{0\})$ (resp.,
$b_1{H_{\Omega_2}}(\mathbf{m};\mathbf{q})\cup (b\Omega_2\times
\{0\})$). By investigating the subset
$$ b_1{H_{\Omega}}(\mathbf{n};\mathbf{p})\cup (b\Omega\times
\{0\})$$ of the boundary $b{H_{\Omega}}(\mathbf{n};\mathbf{p})$ of a
Hua domain $H_{\Omega}(\mathbf{n};\mathbf{p})$, we have that $
b_1{H_{\Omega}}(\mathbf{n};\mathbf{p})\cup (b\Omega\times \{0\})$
consists of $r-\delta$ components
$bPr_j({H_{\Omega}}(\mathbf{n};\mathbf{p}))$ for ${1+\delta}\leq
j\leq r$ (see $(4)$ for the notation of $\delta$) and
$b\Omega\times\{0\}$ is the intersection of these $r-\delta$
components, where
$$Pr_j({H_{\Omega}}(\mathbf{n};\mathbf{p})):={H_{\Omega}}(\mathbf{n};\mathbf{p})\cap
\{w_{(j)}=0\}$$ for $1 + \delta\leq j\leq r.$ Since $f$ is a
biholomorphism between their closures, $f$ maps the subset
$b_1{H_{\Omega_1}}(\mathbf{n};\mathbf{p})\cup (b\Omega_1\times
\{0\})$ of $bH_{\Omega_1}(\mathbf{n};\mathbf{p})$ onto the subset $
b_1{H_{\Omega_2}}(\mathbf{m};\mathbf{q})\cup (b\Omega_2\times
\{0\})$ of $bH_{\Omega_2}(\mathbf{m};\mathbf{q})$. Apply this fact
to $f,\; Pr_{1+\delta}\circ f,\; Pr_{2+\delta}\circ
Pr_{1+\delta}\circ f,\;\cdots,\; Pr_r\circ\cdots\circ
Pr_{1+\delta}\circ f$ in succession, we get
$$f(b\Omega_1\times\{0\})\subset b\Omega_2\times \{0\}.$$
Thus $$f(\Omega_1\times\{0\})\subset\Omega_2\times\{0\}$$ by the
maximum modulus principle. In particular, we have
$f(0,0)\in\Omega_2\times \{0\}$, thus, using fact that $\Gamma
(H_{\Omega_2}(\mathbf{m};\mathbf{q}))$ is transitive on
$\Omega_2\times \{0\}(\subset H_{\Omega_2}(\mathbf{m};\mathbf{q}))$,
we can choose an automorphism $\Phi\in
\Gamma(H_{\Omega_2}(\mathbf{m};\mathbf{q}))$ (see (3)) with
$\Phi(f(0,0))=(0,0)$. Thus $$\Phi\circ
f:H_{\Omega_1}(\mathbf{n};\mathbf{p})\rightarrow
H_{\Omega_2}(\mathbf{m};\mathbf{q})$$ is a biholomorphism with
$\Phi\circ f(0,0)=(0,0)$, therefore, a holomorphic linear
isomorphism by the Cartan's theorem.

At last, we prove that, after a permutation of coordinates, the
$(r+1)\times (r+1)$ block matrix of the holomorphic  linear
isomorphism $\Phi\circ
f:H_{\Omega_1}(\mathbf{n};\mathbf{p})\rightarrow
H_{\Omega_2}(\mathbf{m};\mathbf{q})$ is a block diagonal matrix.
That is, we prove that there exists one and only one nonzero block
in every row of the block matrix. Denote the projection by
$$Pr:H_{\Omega_2}(\mathbf{m};\mathbf{q})\rightarrow \{0\}\times
\Sigma(\mathbf{m};\mathbf{q})
\;(:=H_{\Omega_2}(\mathbf{m};\mathbf{q})\cap
(\{0\}\times\mathbb{C}^{|\mathbf{m}|})).$$ Then we prove that
$$Pr\circ\Phi\circ f\mid_{\overline{\{0\}\times
\Sigma(\mathbf{n};\mathbf{q})}}: \{0\}\times
\Sigma(\mathbf{n};\mathbf{p})\rightarrow \{0\}\times
\Sigma(\mathbf{m};\mathbf{q})$$ must be a holomorphic linear
isomorphism between two generalized complex ellipsoids
$\Sigma(\mathbf{n};\mathbf{p})$ and $\Sigma(\mathbf{m};\mathbf{q})$,
and its matrix $D$ can be obtained from the block matrix of
$\Phi\circ f$ by deleting the first row and first column. In order
to show that $D$ is a block diagonal matrix, we argue by the
contradiction. If there exist no nonzero block or at least two
nonzero blocks $D_{i_1j},D_{i_2j}$ of some column $D_j$ of $D$, then
we have that some strongly pseudoconvex points on
$b\Sigma(\mathbf{n};\mathbf{p})$ are mapped by $Pr\circ\Phi\circ
f\mid_{\overline{\{0\}\times \Sigma(\mathbf{n};\mathbf{q})}}$ to
weakly pseudoconvex points on $b\Sigma(\mathbf{m};\mathbf{q})$. This
is impossible since it is a holomorphic linear isomorphism. Thus,
there exists one and only one nonzero block in every row of the
block matrix $D$. Further, we prove every block except the first one
on the first row and the first column of the matrix of $\Phi\circ f$
is zero. Thus, we get that after a permutation of coordinates, the
$(r+1)\times (r+1)$ block matrix of the linear isomorphism
$\Phi\circ f$ is a block diagonal matrix. These are the key ideas in
proving our main results in this paper.

\section{Preliminaries }
\subsection{Holomorphic extensions of proper holomorphic mappings}
\noindent {\bf Proposition 2.1} {\it Let $\Omega$ be an irreducible
bounded symmetric domain in $\mathbb{C}^d$ of genus $g$ in its
Harish-Chandra realization and let $N_{\Omega}(z,\bar{z})$ be the
generic norm of $\Omega$. Then we have the results as follows:

$(a)$ For any $z_0\in\Omega$, we have
$N_{\Omega}(tz_0,\overline{tz_0})\; (0\leq t\leq 1)$ is a decreasing
function of $t$.

$(b)$ We have
$$N_{\Omega}(z,0)=1 \;\mbox{and} \;\; 0< N_{\Omega}(z,\bar{z})\leq 1\; (z\in \Omega),$$
and $N_{\Omega}(z,\bar{z})=1$ if and only if $z=0$.

$(c)$ Let ${H_{\Omega}}(\mathbf{n};\mathbf{p})$ be a Hua domain.
Then, for any $(z_0,w_{(1)0},\cdots,w_{(r)0})\in
{H_{\Omega}}(\mathbf{n};\mathbf{p})$ and $0\leq t\leq 1$, we have
$(tz_0,tw_{(1)0},\cdots,tw_{(r)0})\in
{H_{\Omega}}(\mathbf{n};\mathbf{p})$. Therefore, each Hua domain is
a starlike domain with respect to the origin of
$\mathbb{C}^{d+|\mathbf{n}|}$, where
$|\mathbf{n}|:=n_1+\cdots+n_r$.}

\vskip 3pt
\noindent
{\bf Proof}. Since $\Omega$ is a bounded circular domain and contains
the origin, there is a homogeneous holomorphic polynomial set
$$\left\{\frac{1}{\sqrt{V(\Omega)}},h_1(z),h_2(z),\cdots \right\},$$
which is an orthonormal basis of the Hilbert space $A^2(\Omega)$ of
square-integrable holomorphic functions on $\Omega$, where
$\deg h_j(z)\geq 1$ (so $h_j(0)=0$) for $j=1,2,\cdots$. Then
\begin{equation}
K_{\Omega}(z,\bar{\xi})=\frac{1}{V(\Omega)}+
h_1(z)\overline{h_1(\xi)}+ h_2(z)\overline{h_2(\xi)}+\cdots
\end{equation}
for all $z, \xi\in \Omega$.

(a) For any $z_0\in \Omega$ and $0\leq t\leq 1$, from (5), we have
$$K_{\Omega}(tz_0,\overline{tz_0})=\frac{1}{V(\Omega)}+t^{2\deg h_1}|h_1(z_0)|^2+t^{2\deg h_2}|h_2(z_0)|^2+\cdots$$
is an increasing function of $t\; (0\leq t\leq 1)$. So
$$N_{\Omega}(tz_0,\overline{tz_0})=(V(\Omega)K_{\Omega}(tz_0,\overline{tz_0}))^{-1/g}$$
is a decreasing function of $t\; (0\leq t\leq 1)$. The proof of
Proposition 2.1 (a) is completed.

(b) Thus, from (5), we have
$$N_{\Omega}(z,0)=(V(\Omega)K_{\Omega}(z,0))^{-1/g}=1$$ for all $z\in
\Omega$ and
$$1\leq V(\Omega)K(z,\bar{z})<+\infty \; (z\in \Omega)\;\mbox{and}\; V(\Omega)K(0,0)=1.$$

Since $\Omega\subset \mathbb{C}^d$ is a bounded circular domain and
contains the origin, we have that
$$\left\{\frac{1}{\sqrt{V(\Omega)}},z_1,(z_2)^2,\cdots,(z_d)^d\right\}$$
is an orthogonal set of the Hilbert space $A^2(\Omega)$ of square-integrable
holomorphic functions on $\Omega$. Take positive numbers $r_1,\cdots,r_d$
such that
$$\left\{\frac{1}{\sqrt{V(\Omega)}},r_1z_1,r_2(z_2)^2,\cdots, r_d(z_d)^d \right\}$$
is an orthonormal set of $A^2(\Omega)$. Then
$$K_{\Omega}(z,\bar{z})\geq \frac{1}{V(\Omega)}+r_1^2|z_1|^2+r_2^2|z_2|^4+\cdots+r_d^2|z_d|^{2d}$$
for all $z=(z_1,\cdots,z_d)\in\Omega$. Thus
$V(\Omega)K_{\Omega}(z,\bar{z})=1\; (z\in\Omega)$ implies $z=0$.

Therefore, we have $0<N_{\Omega}(z,\bar{z})\leq 1 (z\in \Omega)$,
and $N_{\Omega}(z,\bar{z})=1$ if and only if $z=0$. The proof of
Proposition 2.1 (b) is completed.

(c) For any $(z_0,w_{(1)0},\cdots,w_{(r)0})\in
{H_{\Omega}}(\mathbf{n};\mathbf{p})$ and $0\leq t\leq 1$, by
definition, we have
$$\|w_{(1)0}\|^{2p_1}+\cdots+\|w_{(r)0}\|^{2p_r}<N_{\Omega}(z_0,\overline{z_0}).$$
Thus, by (a), we have
$$\|tw_{(1)0}\|^{2p_1}+\cdots+\|tw_{(r)0}\|^{2p_r}\leq \|w_{(1)0}\|^{2p_1}+\cdots+\|w_{(r)0}\|^{2p_r}
<N_{\Omega}(z_0,\overline{z_0})<N_{\Omega}(tz_0,\overline{tz_0}). $$
So we get $(tz_0,tw_{(1)0},\cdots,tw_{(r)0})\in
{H_{\Omega}}(\mathbf{n};\mathbf{p})$. The proof of Proposition 2.1
(c) is completed.

\vskip 10pt
In order to prove that a proper holomorphic mapping between two
equidimensional Hua domains extends holomorphically to their closures,
we need the following lemma.

\vskip 10pt \noindent {\bf Lemma 2.2}.(Bell \cite{Bell}, Theorem 2)
{\it Suppose $f: \Omega_1\rightarrow \Omega_2$ is a proper
holomorphic mapping between bounded circular domains in
$\mathbb{C}^n$. Suppose further that $\Omega_2$ contains the origin
and that the Bergman kernel function $K_{\Omega_1}(z,\bar{\xi})$
associated to $\Omega_1$ is such that for each compact subset $E$ of
$\Omega_1$, there is an open set $U=U(E)$ containing
$\overline{\Omega_1}$ such that $K_{\Omega_1}(z,\bar{\xi})$ extends
to be holomorphic on $U$ as a function of $z$ for each $\xi\in E$.
Then $f$ extends holomorphically to a neighborhood of
$\overline{\Omega_1}$.}

\vskip 10pt Now we prove that a proper holomorphic mapping between
two equidimensional Hua domains extends holomorphically to their
closures as follows (see Lemma 1.1.1 in Mok \cite{Mok} and Th. 2.5
in Tu-Wang \cite{T-W} for references).

\vskip 10pt \noindent {\bf Proposition 2.3}. {\it Let
${H_{\Omega}}(\mathbf{n};\mathbf{p})\subset
\mathbb{C}^{d+|\mathbf{n}|}$ be a Hua domain and $G\subset
\mathbb{C}^{d+|\mathbf{n}|}$ be a bounded circular domain containing
the origin. Suppose that
$F:{H_{\Omega}}(\mathbf{n};\mathbf{p})\rightarrow G$ is a proper
holomorphic mapping. Then $F$ extends holomorphically to a
neighborhood of $\overline{{H_{\Omega}}(\mathbf{n};\mathbf{p})}$.}

\vskip 3pt \noindent {\bf Proof}. Let $r$ be a real number with
$0<r<1$. Since ${H_{\Omega}}(\mathbf{n};\mathbf{p})\subset
\mathbb{C}^{d+|\mathbf{n}|}$ is a starlike domain by Proposition 2.1
(c), we have $rH_{\Omega}(\mathbf{n};\mathbf{p})\subset
{H_{\Omega}}(\mathbf{n};\mathbf{p})$.

Consider the Taylor expansion of the Bergman kernel
$K_{{H_{\Omega}}(\mathbf{n};\mathbf{p})}(z,\bar{\xi})$ on
${H_{\Omega}}(\mathbf{n};\mathbf{p})$ in
$z=(z_1,\cdots,z_{d+|\mathbf{n}|})$ and
$\bar{\xi}=(\bar{\xi_1},\cdots,\overline{\xi_{d+|\mathbf{n}|}})$.
From the invariance of ${H_{\Omega}}(\mathbf{n};\mathbf{p})$ under
the circle group action $z\rightarrow
e^{\sqrt{-1}\theta}z\;(\theta\in\mathbb{R})$, we have the Bergman
kernel $K_{{H_{\Omega}}(\mathbf{n};\mathbf{p})}(z,\bar{\xi})$ on
${H_{\Omega}}(\mathbf{n};\mathbf{p})$ is invariant under the circle
group action. It follows that the coefficient of
$z^I\overline{\xi^J}$ is zero whenever $|I|\neq |J|$. Thus, the
Bergman kernel
$K_{{H_{\Omega}}(\mathbf{n};\mathbf{p})}(z,\bar{\xi})$ on
${H_{\Omega}}(\mathbf{n};\mathbf{p})$ is of the form
$$K_{{H_{\Omega}}(\mathbf{n};\mathbf{p})}(z,\bar{\xi})=\sum_{|I|=|J|}a_{I\bar{J}}z^I\overline{\xi^J}$$
for $z,\xi\in {H_{\Omega}}(\mathbf{n};\mathbf{p})$. Since
$(rz)^I\overline{(\xi/r)^J}=z^I\overline{\xi^J}$ whenever $|I|=|J|$,
we have
$$K_{{H_{\Omega}}(\mathbf{n};\mathbf{p})}(z,\bar{\xi})=K_{{H_{\Omega}}(\mathbf{n};\mathbf{p})}(rz,\overline{\xi/r})$$
for all $z\in {H_{\Omega}}(\mathbf{n};\mathbf{p})$, $\xi\in
rH_{\Omega}(\mathbf{n};\mathbf{p})$. Then, for every fixed $\xi\in
rH_{\Omega}(\mathbf{n};\mathbf{p})$, we have
$K_{{H_{\Omega}}(\mathbf{n};\mathbf{p})}(z,\bar{\xi})$ extends
holomorphically to $\frac{1}{r}{H_{\Omega}}(\mathbf{n};\mathbf{p})$
as a function of $z$.

Therefore, for each compact subset $E$ of
${H_{\Omega}}(\mathbf{n};\mathbf{p})$, there exists a real number
$r_0\; (0<r_0<1)$ with $E \subset r_0
{H_{\Omega}}(\mathbf{n};\mathbf{p})$ such that
$K_{{H_{\Omega}}(\mathbf{n};\mathbf{p})}(z,\bar{\xi})$ extends
holomorphically to $\frac{1}{r_0}
{H_{\Omega}}(\mathbf{n};\mathbf{p})$ (a neighborhood of
$\overline{{H_{\Omega}}(\mathbf{n};\mathbf{p})}$) as a function of
$z$ for all $\xi\in E$. By Lemma 2.2, we have that $f$ extends
holomorphically to a neighborhood of
$\overline{{H_{\Omega}}(\mathbf{n};\mathbf{p})}$. The proof of
Proposition 2.3 is finished.

\subsection{The structure of the boundary of a Hua domain ${H_{\Omega}}(\mathbf{n},\mathbf{p})$}

For a Hua domain
${H_{\Omega}}(\mathbf{n};\mathbf{p})={H_{\Omega}}(n_1,\cdots,n_r;p_1,\cdots,p_r)$
in its standard form, we will investigate the strongly pseudoconvex
part of its boundary $bH_{\Omega}(\mathbf{n};\mathbf{p})$ which is
comprised of
\begin{align*}
&bH_{\Omega}(\mathbf{n};\mathbf{p})
=b_0{H_{\Omega}}(\mathbf{n};\mathbf{p})\cup
b_1{H_{\Omega}}(\mathbf{n};\mathbf{p})\cup (b\Omega\times \{0\}),
\end{align*}
where $b_0{H_{\Omega}}(\mathbf{n};\mathbf{p})$ and
$b_1{H_{\Omega}}(\mathbf{n};\mathbf{p})$ are the same as those in
\eqref{1.4}.

\vskip 10pt \noindent {\bf Proposition 2.4}. {\it Let $\Omega\subset
\mathbb{C}^d$ be an irreducible bounded symmetric domain of genus
$g$ in its Harish-Chandra realization. Then we have the conclusions
as follows.

$(a)$ $b_0{H_{\Omega}}(\mathbf{n};\mathbf{p})$ is a real analytic
hypersurface in $\mathbb{C}^{d+|\mathbf{n}|}$ and
${H_{\Omega}}(\mathbf{n};\mathbf{p})$ is strongly pseudoconvex at
all points of $b_0{H_{\Omega}}(\mathbf{n};\mathbf{p})$.

$(b)$ If ${H_{\Omega}}(\mathbf{n};\mathbf{p})$ isn't a ball, then
${H_{\Omega}}(\mathbf{n};\mathbf{p})$ is not strongly pseudoconvex
at any point of $$b_1{H_{\Omega}}(\mathbf{n};\mathbf{p})\cup
(b\Omega\times \{0\}).$$ Obviously,
$b_1{H_{\Omega}}(\mathbf{n};\mathbf{p})\cup (b\Omega\times \{0\})$
is  contained in a complex analytic set of complex codimension
$\min\{n_{1+\delta},\cdots,n_r,n_{1}+\cdots+n_r\}$. (cf. Lemma 3.2
in Rong \cite{R}.) }

\vskip 10pt \noindent {\bf Proof}. Let $\{h_i(z)\}_{i=1}^{\infty}$
be an orthonormal basis of the Hilbert space $A^2(\Omega)$ of
square-integrable holomorphic functions. Then we have
$$K_{\Omega}(z,\bar{z})=\sum_{i=1}^{\infty}h_i(z)\overline{h_i(z)}$$
converges uniformly on any compact subset of $\Omega$. Let
$$\rho(z,w_{(1)},\cdots,w_{(r)}):=\|w_{(1)}\|^{2p_1}+\cdots+\|w_{(r)}\|^{2p_r}-\sigma
(K_{\Omega}(z,\bar{z}))^{-\lambda}$$ where
$\sigma:=(V(\Omega))^{-\frac{1}{g}}$ and $\lambda:=\frac{1}{g}$ are
positive. Then $\rho$ is a real analytic definition function of
$b_0{H_{\Omega}}(\mathbf{n};\mathbf{p})$.

Fix a point $(z_0,w_{(1)0},\cdots,w_{(r)0})\in
b_0{H_{\Omega}}(\mathbf{n};\mathbf{p}) (\subset
\mathbb{C}^{d}\times\mathbb{C}^{n_1}\times\cdots\times\mathbb{C}^{n_r})$
and let
$$T=(\xi,\eta_1,\cdots,\eta_r)\in T^{1,0}_{(z_0,w_{(1)0},\cdots,w_{(r)0})}
(b_0{H_{\Omega}}(\mathbf{n};\mathbf{p}))(\subset
\mathbb{C}^{d}\times\mathbb{C}^{n_1}\times\cdots\times\mathbb{C}^{n_r}).$$
Then by definition, we have
\begin{align}
&w_{(j)0}\neq 0,\;\; j=1,\cdots,r;\\
&\|w_{(1)0}\|^{2p_1}+\cdots+\|w_{(r)0}\|^{2p_r}-\sigma (K_{\Omega}(z_0,\bar{z_0}))^{-\lambda}=0;\\
&\sum_{k=1}^r p_k\|w_{(k)0}\|^{2(p_k-1)}(\overline{w_{(k)0}}\cdot \eta_k)
+\sigma\lambda(K_{\Omega}(z_0,\bar{z_0}))^{-\lambda-1}\sum_{i=1}^{\infty}\overline{h_i(z_0)}(h'_i(z_0)\cdot\xi)=0,
\end{align}
where $h'_i(z_0)\cdot \xi=\sum_{k=1}^d\frac{\partial h_i}{\partial z_k}(z_0)\xi_k$.

Therefore, from (6),(7),(8), the Levi form of $\rho$ at the point
$(z_0,w_{(1)0},\cdots,w_{(r)0})$ is computed as follows:
\begin{align*}
L_{\rho}&(T,T):=\sum_{i,j=1}^{d+|\mathbf{n}|}
\frac{\partial^2 \rho}{\partial T_i\partial\overline{T_j}}(z_0,w_{(1)0},\cdots,w_{(r)0})T_i\overline{T_j}\\
=&\sum_{k=1}^r p_k\|w_{(k)0}\|^{2(p_k-1)}\|\eta_k\|^2
 +\sum_{k=1}^r p_k(p_k-1)\|w_{(k)0}\|^{2(p_k-2)}|\overline{w_{(k)0}}\cdot\eta_k|^2\\
&+\sigma\lambda K_{\Omega}(z_0,\overline{z_0})^{-(\lambda+2)}
[K_{\Omega}(z_0,\overline{z_0})\sum_{i=1}^{\infty}|h'_i(z_0)\cdot\xi|^2-(\lambda+1)|
\sum_{i=1}^{\infty}\overline{h_i(z_0)}(h'_i(z_0)\cdot\xi)|^2]\\
=&\sum_{k=1}^r p_k(\|w_{(k)0}\|^{2(p_k-1)}\|\eta_k\|^2 -\|w_{(k)0}\|^{2(p_k-2)}|\overline{w_{(k)0}}\cdot\eta_k|^2)\\
&+\sum_{k=1}^r p_k^2\|w_{(k)0}\|^{2(p_k-2)}|\overline{w_{(k)0}}\cdot\eta_k|^2
 -\sigma\lambda^2 K_{\Omega}(z_0,\overline{z_0})^{-(\lambda+2)} |\sum_{i=1}^{\infty}\overline{h_i(z_0)}(h'_i(z_0)\cdot\xi)|^2 \\
 &+\sigma\lambda K_{\Omega}(z_0,\overline{z_0})^{-(\lambda+2)}
\Big[K_{\Omega}(z_0,\overline{z_0})\sum_{i=1}^{\infty}|h'_i(z_0)\cdot\xi|^2-|\sum_{i=1}^{\infty}\overline{h_i(z_0)}
(h'_i(z_0)\cdot\xi)|^2\Big]
\end{align*}
\begin{align*}
=&\sum_{k=1}^r p_k\|w_{(k)0}\|^{2(p_k-2)}\Big[\|w_{(k)0}\|^2\|\eta_k\|^2-|\overline{w_{(k)0}}\cdot\eta_k|^2\Big]+(\sum_{k=1}^r\|w_{(k)0}\|^2)^{-1}\\
&\times \Big[(\sum_{k=1}^r p_k^2\|w_{(k)0}\|^{2(p_k-2)}|\overline{w_{(k)0}}\cdot\eta_k|^2)(\sum_{k=1}^r\|w_{(k)0}\|^2)
-\big|\sum_{k=1}^r p_k \|w_{(k)0}\|^{2(p_k-1)}(\overline{w_{(k)0}}\cdot \eta_k)\big|^2\Big]\\
&+\sigma\lambda K_{\Omega}(z_0,\overline{z_0})^{-(\lambda+2)}
\Big[(\sum_{i=1}^{\infty}|h_i(z_0)|^2)(\sum_{i=1}^{\infty}|h'_i(z_0)\cdot\xi|^2)-|\sum_{i=1}^{\infty}\overline{h_i(z_0)}(h'_i(z_0)\cdot\xi)|^2\Big]\\
\geq& 0
\end{align*}
by the Cauchy-Schwarz inequality, for all
$T=(\xi,\eta_1,\cdots,\eta_r)\in
T^{1,0}_{(z_0,w_{(1)0},\cdots,w_{(r)0})}(b_0{H_{\Omega}}(\mathbf{n};\mathbf{p}))$
and the equality holds if and only if
\begin{align}
& \|w_{(k)0}\|^2\|\eta_k\|^2-|\overline{w_{(k)0}}\cdot\eta_k|^2=0 ,\;\;\; 1\leq k\leq r;\\
& (\sum_{k=1}^r p_k^2\|w_{(k)0}\|^{2(p_k-2)}|\overline{w_{(k)0}}\cdot\eta_k|^2)(\sum_{k=1}^r\|w_{(k)0}\|^2)
-\big|\sum_{k=1}^r p_k \|w_{(k)0}\|^{2(p_k-1)}(\overline{w_{(k)0}}\cdot \eta_k)\big|^2=0;\\
& (\sum_{i=1}^{\infty}|h_i(z_0)|^2)(\sum_{i=1}^{\infty}|h'_i(z_0)\cdot\xi|^2)-|\sum_{i=1}^{\infty}\overline{h_i(z_0)}(h'_i(z_0)\cdot\xi)|^2=0.
\end{align}

Now we prove the Levi form $L_{\rho}(T,T)$ of $\rho$ at the point
$(z_0,w_{(1)0},\cdots,w_{(r)0})$ is positive for any nonzero
$T=(\xi,\eta_1,\cdots,\eta_r)\in
T^{1,0}_{(z_0,w_{(1)0},\cdots,w_{(r)0})}(b_0{H_{\Omega}}(\mathbf{n};\mathbf{p}))$
as follows:

{\bf Case 1} Suppose $\xi\neq 0$. Since
$$K_{\Omega}(z,\bar{z})=\sum_{i=1}^{\infty}h_i(z)\overline{h_i(z)}$$
is independent of the choice of the orthonormal basis
$\{h_i(z)\}_{i=1}^{\infty}$ of the Hilbert space $A^2(\Omega)$ and
$\Omega$ is bounded, we may choose that $h_1(z)$ is a nonzero
constant and $h_2(z)$ satisfies $h'_2(z_0)\cdot \xi\neq 0$. This
gives that $(h_1(z_0),h_2(z_0))$ and $(h_1'(z_0)\cdot \xi,
h_2'(z_0)\cdot \xi)$ are linearly independent. Thus
$$\sum_{k=1}^{\infty}|h_k(z_0)|^2\sum_{i=1}^{\infty}|h'_i(z_0)\cdot\xi|^2-
|\sum_{i=1}^{\infty}\overline{h_i(z_0)}(h'_i(z_0)\cdot\xi)|^2>0$$
which is a contradiction with (11). Therefore, $L_{\rho}(T,T)>0$ for
all $$T=(\xi,\eta_1,\cdots,\eta_r)\in
T^{1,0}_{(z_0,w_{(1)0},\cdots,w_{(r)0})}\\(b_0{H_{\Omega}}(\mathbf{n};\mathbf{p}))$$
with $\xi\neq 0$.

{\bf Case 2} Suppose $\xi=0$. Then $T=(\xi,\eta_1,\cdots,\eta_r)\neq
0$ implies that there exits $\eta_{i_0}\neq 0$. On the other hand,
since $\xi=0$, by (8), we have
$$\sum_{k=1}^r p_k\|w_{(k)0}\|^{2(p_k-1)}(\overline{w_{(k)0}}\cdot \eta_k)=0. $$
Hence, by (10), we get
$$(\sum_{k=1}^r p_k^2\|w_{(k)0}\|^{2(p_k-1)}\|\eta_k\|^2)(\sum_{k=1}^r\|w_{(k)0}\|^2)=0.$$
Since $\|w_{(j)0}\|^{2}\neq 0,1\leq j\leq r$, we have $\eta_{i}=0,
1\leq i\leq r$, this is a contradiction. Therefore,
$L_{\rho}(T,T)>0$ for all $$T=(0,\eta_1,\cdots,\eta_r)\in
T^{1,0}_{(z_0,w_{(1)0},\cdots,w_{(r)0})}(b_0{H_{\Omega}}(\mathbf{n};\mathbf{p}))$$
with $(\eta_1,\cdots, \eta_r)\neq (0,\cdots,0)$.

Thus the Levi form $L_{\rho}(T,T)$ is positive definite on
$T^{1,0}_{(z_0,w_{(1)0},\cdots,w_{(r)0})}(b_0{H_{\Omega}}(\mathbf{n};\mathbf{p}))$.
This means that every point of
$b_0{H_{\Omega}}(\mathbf{n};\mathbf{p})$ is strongly pseudoconvex.

Let $(z_0,w_{(1)0},\cdots,w_{(r)0})\in
b_1{H_{\Omega}}(\mathbf{n};\mathbf{p})$. Without loss of generality,
we assume  $\|w_{(k)0}\|=0$ for $1+\delta\leq k\leq i_0$ and
$\|w_{(k)0}\|\not=0$ for $i_0+1\leq k\leq r$. Suppose
 $p_k\geq 2$ for $1+\delta\leq k\leq i_0$ (otherwise, $b{H_{\Omega}}(\mathbf{n};\mathbf{p})$ is not $C^2$
 at the point $(z_0,w_{(1)0},\cdots,w_{(r)0})$, so $b{H_{\Omega}}(\mathbf{n};\mathbf{p})$ is not
  strongly pseudoconvex at the point). In this case, take $T_0=(0,0,\cdots,\eta_{i_0},
\cdots,0)\in
T^{1,0}_{(z_0,w_{(1)0},\cdots,w_{(r)0})}(bH_{\Omega}(\mathbf{n};\mathbf{p}))$
with $\|\eta_{i_0}\|\neq 0$, then $L_{\rho}(T_0, T_0)=0$. Hence,
${H_{\Omega}}(\mathbf{n};\mathbf{p})$ is not strongly pseudoconvex
at any point of $b_1{H_{\Omega}}(\mathbf{n};\mathbf{p})$.

For any irreducible bounded symmetric domain $\Omega$ in
$\mathbb{C}^d$, we have
$\Gamma({H_{\Omega}}(\mathbf{n};\mathbf{p}))$ is transitive on
$\Omega\times\{0\}(\subset {H_{\Omega}}(\mathbf{n};\mathbf{p}))$.
Since ${H_{\Omega}}(\mathbf{n};\mathbf{p})$ is not the unit ball,
${H_{\Omega}}(\mathbf{n};\mathbf{p})$ is not strictly pseudoconvex
at any point of $b\Omega\times\{0\}$ by the Wong-Rosay theorem. The
proof of Proposition 2.4 is completed.

\vskip 10pt \noindent {\bf Lemma 2.5}. (Pinchuk \cite{Pin}, Lemma
1.3) {\it Let $D_1,D_2\subset \mathbb{C}^n$ be two domains, $p\in
bD_1$, and let $U$ be a neighborhood of $p$ in $\mathbb{C}^n$ such
that $U\cap \overline{D_1}$ is connected. Suppose that the mapping
$f=(f_1,\cdots,f_n):U\cap \overline{D_1}\rightarrow \mathbb{C}^n$ is
continuously differentiable in $U\cap \overline{D_1}$ and
holomorphic in $U\cap D_1$ with $f(U\cap bD_1)\subset bD_2$. Take a
domain $V\subset \mathbb{C}^n$ with $f(U\cap D_1)\subset V$. Suppose
that $U\cap bD_1$ and $U\cap bD_2$ are strongly pseudoconvex
hypersurfaces in $\mathbb{C}^n$. Then either $f$ is constant or the
holomorphic Jacobian  determinant $J_f(z)=\det(\frac{\partial
f_i}{\partial z_j})$ does not vanish in $U\cap bD_1$.}

\vskip 10pt \noindent {\bf Lemma 2.6}. {\it Let Hua domains
$H_{\Omega_1}(\mathbf{n};\mathbf{p})$ and
$H_{\Omega_2}(\mathbf{m};\mathbf{q})$ be in their standard forms,
where $\Omega_1\subset \mathbb{C}^{d_1}$ and $\Omega_2\subset
\mathbb{C}^{d_2}$ are two irreducible bounded symmetric domains in
the Harish-Chandra realization, and $\mathbf{n},
\mathbf{m}\in\mathbb{N}^r$,
$\mathbf{p},\mathbf{q}\in(\mathbb{R}_{+})^r$.  Then every
biholomorphism $f:H_{\Omega_1}(\mathbf{n},\mathbf{p})\rightarrow
H_{\Omega_2}(\mathbf{m},\mathbf{q})$ sends $\Omega_1\times\{0\}$
into $\Omega_2\times\{0\}$. Therefore, we have
$f(\Omega_1\times\{0\}) =\Omega_2\times\{0\}$. }

\vskip 3pt \noindent {\bf Proof}. We will divide our proof into the
following two cases.

Case 1. Suppose that $H_{\Omega_1}(\mathbf{n},\mathbf{p})$ is a unit
ball. Since there exists a biholomorphism
$f:H_{\Omega_1}(\mathbf{n},\mathbf{p})\rightarrow
H_{\Omega_2}(\mathbf{m},\mathbf{q})$, we have
$H_{\Omega_2}(\mathbf{m},\mathbf{q})$ must be the unit ball also.

In fact, since there exists a biholomorphism $f$ from the unit ball
onto $H_{\Omega_2}(\mathbf{m},\mathbf{q})$, we have
$H_{\Omega_2}(\mathbf{m},\mathbf{q})$ must be a bounded symmetric
domain with rank 1.  Then  $b\Omega_2 \times\{0\} \;(\subset
bH_{\Omega_2}(\mathbf{m},\mathbf{q}))$ can not contain any
positive-dimensional complex submanifold, and so $\Omega_2$ is a
bounded symmetric domain with rank 1 in the Harish-Chandra
realization. This means $\Omega_2$ is a unit ball. So we have
$${H_{\Omega_2}}(\mathbf{m};\mathbf{q})=\left\{(z,w_{(1)},\cdots,w_{(r)})\in B^d \times
\mathbb{C}^{m_1}\times\cdots\times\mathbb{C}^{m_r}: \|z\|^2+
\sum_{j=1}^r\|w_{(j)}\|^{2q_j}< 1\right\}.$$ It is known that the
generalized complex ellipsoid is homogeneous if and only if $q_j=1$
for all $1\leq j\leq r$ (cf. Kodama \cite{K1}). Thus we have
$H_{\Omega_2}(\mathbf{m},\mathbf{q})$ must be the unit ball.

From Hua domains $H_{\Omega_1}(\mathbf{n};\mathbf{p})$ and
$H_{\Omega_2}(\mathbf{m};\mathbf{q})$ being in their standard forms,
we get $H_{\Omega_1}(\mathbf{n},\mathbf{p})= \Omega_1\;(\cong
\Omega_1\times\{0\})$ and
$H_{\Omega_2}(\mathbf{m},\mathbf{q})=\Omega_2 \;(\cong
\Omega_2\times\{0\}).$ Then, Lemma 2.6 is true.

Case 2. Suppose that $H_{\Omega_1}(\mathbf{n},\mathbf{p})$ is not a
unit ball. Let $f:{H_{\Omega_1}}(\mathbf{n};\mathbf{p})\rightarrow
{H_{\Omega_2}}(\mathbf{m};\mathbf{q})$ be a biholomorphism. By
Proposition 2.3, the biholomorphism
$$f:H_{\Omega_1}(\mathbf{n};\mathbf{p})\rightarrow
H_{\Omega_2}(\mathbf{m};\mathbf{q})$$ extends a biholomorphism
between $\overline{H_{\Omega_1}(\mathbf{n};\mathbf{p})}$ and
$\overline{H_{\Omega_2}(\mathbf{m};\mathbf{q})}$ by the uniqueness
theorem. So, by Proposition 2.4(b), we have
$$f((b\Omega_1\times \{0\})\cup
b_1H_{\Omega_1}(\mathbf{n};\mathbf{p}))=
(b\Omega_2\times \{0\})\cup
b_1H_{\Omega_2}(\mathbf{m};\mathbf{q}).$$
Since
$$(b\Omega_1\times\{0\})\cup b_1H_{\Omega_1}(\mathbf{n};\mathbf{p})
=\bigcup_{j=1+\delta}^r bPr_j(H_{\Omega_1}(\mathbf{n};\mathbf{p}))$$
and $$(b\Omega_2\times \{0\})\cup
b_1H_{\Omega_2}(\mathbf{m};\mathbf{q})=\bigcup_{j=1+\varepsilon}^r b
Pr_j(H_{\Omega_2}(\mathbf{m};\mathbf{q})),$$ where
$$Pr_j(H_{\Omega_1}(\mathbf{n};\mathbf{p})):=H_{\Omega_1}(\mathbf{n};\mathbf{p})\cap
\{w_{(j)}=0\}\;\;(1+\delta\leq j\leq r)$$ and
$$Pr_j(H_{\Omega_2}(\mathbf{m};\mathbf{q})):=H_{\Omega_2}(\mathbf{m};\mathbf{q})\cap
\{w'_{(j)}=0\}\;\;(1+\varepsilon\leq j\leq r),$$ in which
\begin{equation*}
\delta=
\begin{cases}
1 & \text{if $p_1=1$}, \\
0 & \text{if $p_1\neq 1$},
\end{cases}\;\;\;
\varepsilon=
\begin{cases}
1 & \text{if $q_1=1$}, \\
0 & \text{if $q_1\neq 1$},
\end{cases}
\end{equation*}
we have
$$f(\bigcup_{j=1+\delta}^r b
Pr_j(H_{\Omega_1}(\mathbf{n};\mathbf{p})))=
\bigcup_{j=1+\varepsilon}^r b
Pr_j(H_{\Omega_2}(\mathbf{m};\mathbf{q})).$$
In particular, we have
$$f(b Pr_{1+\delta}(H_{\Omega_1}(\mathbf{n};\mathbf{p})))\subset
\bigcup_{j=1+\varepsilon}^r b
Pr_j(H_{\Omega_2}(\mathbf{m};\mathbf{q})).$$

Let $f=(\tilde{f},f_{(1)},\cdots,f_{(r)})$ and $U_j:=\{\xi\in
\overline{Pr_{1+\delta}(H_{\Omega_1}(\mathbf{n},\mathbf{p}))}:
f_{(j)}(\xi)=0\}.$  Then $$b
Pr_{1+\delta}(H_{\Omega_1}(\mathbf{n};\mathbf{p}))\subset
\bigcup_{j=1+\epsilon}^r U_j.$$ Since $b
Pr_{1+\delta}(H_{\Omega_1}(\mathbf{n};\mathbf{p}))$ has real
codimension $1$, there exists a $U_{j_0}$ with real codimension
$\leq 1$. On the other hand, if $U_{j_0}$ is a proper complex
analytic subset of
$\overline{Pr_{1+\delta}(H_{\Omega_1}(\mathbf{n},\mathbf{p}))}$,
then $U_{j_0}$ has real codimension $\geq 2$. Thus there exists a
$U_{j_0}$ such that
$$U_{j_0}=\overline{Pr_{1+\delta}(H_{\Omega_1}(\mathbf{n},\mathbf{p}))}.$$
That is, $f_{(j_0)}\equiv 0$ on
$Pr_{1+\delta}(H_{\Omega_1}(\mathbf{n};\mathbf{p}))$. Hence,
\begin{align*}
f\mid_{Pr_{1+\delta}(H_{\Omega_1}(\mathbf{n};\mathbf{p}))}:Pr_{1+\delta}(H_{\Omega_1}(\mathbf{n};\mathbf{p}))\rightarrow Pr_{j_0}(H_{\Omega_2}
(\mathbf{m};\mathbf{q}))
\end{align*}
is a proper holomorphic mapping and holomorphic on the closure of
$Pr_{1+\delta}(H_{\Omega_1}(\mathbf{n};\mathbf{p})).$ If
$n_{1+\delta}<m_{j_0}$, then $
\dim{Pr_{1+\delta}(H_{\Omega_1}(\mathbf{n};\mathbf{p}))} >
\dim{Pr_{j_0}(H_{\Omega_2}(\mathbf{m};\mathbf{q}))}$. Since
$f\mid_{Pr_{1+\delta}(H_{\Omega_1}(\mathbf{n};\mathbf{p}))}$ is
proper, it is
 a contradiction. Thus, $n_{1+\delta}\geq m_{j_0}$.

Since $f$ is a biholomorphism, by the similar argument, we have
$f^{-1}( Pr_{j_0}(H_{\Omega_2} (\mathbf{m};\mathbf{q}))) \subset
Pr_{1+\delta}(H_{\Omega_1}(\mathbf{n};\mathbf{p}))$ and
\begin{align*}
f^{-1}\mid_{Pr_{j_0}(H_{\Omega_2} (\mathbf{m};\mathbf{q}))}:
Pr_{j_0}(H_{\Omega_2} (\mathbf{m};\mathbf{q})) \rightarrow
Pr_{1+\delta}(H_{\Omega_1}(\mathbf{n};\mathbf{p}))
\end{align*}
is a proper holomorphic mapping and holomorphic on the closure of
$Pr_{j_0}(H_{\Omega_2} (\mathbf{m};\mathbf{q})).$ Therefore, we have
$m_{j_0} \geq n_{1+\delta}$ and so we have  $n_{1+\delta}= m_{j_0}$.

This means that
\begin{align*}
f\mid_{Pr_{1+\delta}(H_{\Omega_1}(\mathbf{n};\mathbf{p}))}:Pr_{1+\delta}(H_{\Omega_1}(\mathbf{n};\mathbf{p}))\rightarrow
Pr_{j_0}(H_{\Omega_2} (\mathbf{m};\mathbf{q}))
\end{align*}
is a biholomorphism and holomorphic on the closure of
$Pr_{1+\delta}(H_{\Omega_1}(\mathbf{n};\mathbf{p})).$ Thus, we have
$$f\mid_{\overline{Pr_{1+\delta}(H_{\Omega_1}(\mathbf{n};\mathbf{p}))}}\big((b\Omega_1\times\{0\})\cup b_1Pr_{1+\delta}
(H_{\Omega_1}(\mathbf{m},\mathbf{q}))\big)
\subset\big(b\Omega_2\times\{0\})\cup
b_1Pr_{j_0}(H_{\Omega_2}(\mathbf{m},\mathbf{q})\big)$$ by
Proposition 2.4(b). Therefore,
\begin{align*}
f\mid_{\overline{Pr_{1+\delta}(H_{\Omega_1}(\mathbf{n};\mathbf{p}))}}
(\bigcup_{j=2+\delta}^{r}
bPr_j(Pr_1(H_{\Omega_1}(\mathbf{n};\mathbf{p}))))\subset
\bigcup_{j=1+\varepsilon,j\neq j_0}^{r}
bPr_j(Pr_{j_0}(H_{\Omega_2}(\mathbf{m};\mathbf{q}))),
\end{align*}
where $$Pr_j(Pr_1(H_{\Omega_1}(\mathbf{n};\mathbf{p})))=
Pr_1(H_{\Omega_1}(\mathbf{n};\mathbf{p}))\cap \{w_{(j)}=0\}\;\;
(2+\delta\leq j\leq r).$$

By induction, since $b\Omega_1\times \{0\}=Pr_{r}(\cdots
Pr_{2+\delta}(Pr_{1+\delta}H_{\Omega_1}(\mathbf{n};\mathbf{p})))$,
we have
$f\mid_{\overline{\Omega_1\times\{0\}}}(b\Omega_1\times\{0\})\subset
(b\Omega_2\times \{0\})$ and thus
$$f(\Omega_1\times\{0\})\subset \Omega_2\times \{0\}$$
by the maximum modulus principle. The proof of Lemma 2.6 is completed.

\vskip 6pt Remark. It is important that the Hua domain is written
{\it in its standard form} in Lemma 2.6. For example, define
\begin{align*}
H_{\mathbf{B}^2}(\mathrm{(2,2)};\mathrm{(1,2)})=\left\{(z,w_{(1)},w_{(2)}
)\in \mathbf{B}^2 \times \mathbb{C}^2 \times  \mathbb{C}^2: \|z\|^2+
\|w_{(1)}\|^{2} + \|w_{(2)}\|^{4}< 1\right\}.
\end{align*}
Then
$H_{\mathbf{B}^2}(\mathrm{(2,2)};\mathrm{(1,2)})=H_{\mathbf{B}^4}(\mathrm{(2)};\mathrm{(2)}),$
and an automorphism
$$\varphi\in {\rm
Aut}(H_{\mathbf{B}^2}(\mathrm{(2,2)};\mathrm{(1,2)}))\;(={\rm
Aut}(H_{\mathbf{B}^4}(\mathrm{(2)};\mathrm{(2)})))
$$
sends $\mathbf{B}^4\times\{0\}$ into $\mathbf{B}^4\times\{0\},$ but
in general, can not send $\mathbf{B}^2\times\{0\}$ into
$\mathbf{B}^2\times\{0\}$ (see Theorem 1.B for references).

\subsection{Complex linear isomorphisms between two equidimensional generalized complex ellipsoids}
In order to get the explicit form of the biholomorphisms  between
two equidimensional Hua domains, we need following two lemmas about
generalized complex ellipsoids.

\vskip 10pt \noindent {\bf Lemma 2.7}. {\it Let $D_j\in M_{n_j\times
m},(1\leq j\leq r)$ be $n_j\times m$ matrix with $m\leq n_j$ such
that the $(n_1+\cdots+n_r)\times m$ matrix $(D_1,\cdots,D_r)^t$ has
rank $m$. If the system of linear equations
\begin{equation}
\sum_{j=1}^r\zeta_jD_j=\mathbf{0},
\end{equation}
where $\zeta_j\in\mathbb{C}^{n_j}$ $(1\leq j\leq r),$ does not have
solution $(\alpha_1,\cdots,\alpha_r)\in
(\mathbb{C}^{n_1}\setminus\{\mathbf{0}\})\times\cdots\times(\mathbb{C}^{n_r}\setminus\{\mathbf{0}\})$,
then there exists at least one $n_j$ such that $n_j=m$ and only one
$D_{j_0}$ with $D_{j_0}\neq \mathbf{0}$. Moreover, $n_{j_0}=m$ and
$D_{j_0}$ is a nonsingular $m\times m$ matrix.}

\vskip 3pt \noindent {\bf Proof}. Suppose that each $n_j>m,
j=1,\cdots,r$, then each system of linear equations
$$\zeta_jD_j=\mathbf{0}$$ has a solution $\alpha_j\in
\mathbb{C}^{n_j}\setminus\{\mathbf{0}\}$. Thus
$(\alpha_1,\cdots,\alpha_r)\in
(\mathbb{C}^{n_1}\setminus\{\mathbf{0}\})\times\cdots\times(\mathbb{C}^{n_r}\setminus\{\mathbf{0}\})$
is a solution of the system of linear equations (12), a
contradiction. Hence there exists at least one $n_j$ such that
$n_j=m$.

If all $D_j's$ are singular $m\times m$ matrices, then by the same
reasoning as above, we can get a solution
$(\alpha_1,\cdots,\alpha_r)\in
(\mathbb{C}^{n_1}\setminus\{\mathbf{0}\})\times\cdots\times(\mathbb{C}^{n_r}\setminus\{\mathbf{0}\})$
of the system of linear equations (12), a contradiction. Thus there
exists a nonsingular $m\times m$ matrix, say $D_1$.

If there exists another $D_j$ with $D_j\neq \mathbf{0}$, then we can
choose $(\alpha_2,\cdots,\alpha_r)\in
(\mathbb{C}^{n_2}\setminus\{\mathbf{0}\})\times\cdots\times(\mathbb{C}^{n_r}\setminus\{\mathbf{0}\})$,
such that
$$\sum_{j=2}^{r}\alpha_jD_j\neq \mathbf{0}.$$
Consider the system of linear equations
$$\zeta_1D_1=\sum_{j=2}^{r}\alpha_jD_j.$$
Since $D_1$ is nonsingular and $ \sum_{j=2}^{r}\alpha_jD_j\neq
\mathbf{0}$, it has a unique solution
$\alpha_1\in\mathbb{C}^{n_1}\setminus\{\mathbf{0}\}$. Thus
$(\alpha_1,\cdots,\alpha_r)\in
(\mathbb{C}^{n_1}\setminus\{\mathbf{0}\})\times\cdots\times(\mathbb{C}^{n_r}\setminus\{\mathbf{0}\})$
is a solution of the system of linear equations (12), a
contradiction. Thus, $D_1$ is the unique nonzero matrix. The proof
of Lemma 2.7 is finished.

\vskip 6pt \noindent {\bf Lemma 2.8}. {\it Let
$\Sigma(\mathbf{n};\mathbf{p})$ and $\Sigma(\mathbf{m};\mathbf{q})$
be two equidimensional generalized complex ellipsoids, where
$\mathbf{n},\mathbf{m}\in\mathbb{N}^r$,
$\mathbf{p},\mathbf{q}\in(\mathbb{R}_{+})^r$ (where $p_k\neq 1,\;
q_k\neq 1$ for $2\leq k\leq r$). Let
$h:\Sigma(\mathbf{n};\mathbf{p})\rightarrow
\Sigma(\mathbf{m};\mathbf{q})$ be a biholomorphic linear isomorphism
between $\Sigma(\mathbf{n};\mathbf{p})$ and
$\Sigma(\mathbf{m};\mathbf{q})$. Then there exists a permutation
$\sigma\in S_r$ such that $n_{\sigma(i)}=m_i, \; p_{\sigma(i)}=q_i$
and
\begin{eqnarray*}
\begin{aligned}
h(\zeta_1,\zeta_2,\cdots,\zeta_r)=(\zeta_{\sigma(1)},\zeta_{\sigma(2)},\cdots,\zeta_{\sigma(r)})
\begin{pmatrix}
U_1 &     &       &    \\
    & U_2 &       &    \\
    &     &\ddots &    \\
    &     &       & U_r\\
\end{pmatrix}
\end{aligned},
\end{eqnarray*}
where $U_i$ is a unitary transformation of
$\mathbb{C}^{m_i}(m_i=n_{\sigma(i)})$ for $1\leq i\leq r$.}

\vskip 6pt \noindent {\bf Proof}.  Let
\begin{equation*}
\delta=
\begin{cases}
1 & \text{if $p_1=1$} \\
0 & \text{if $p_1\neq 1$}
\end{cases},\;\;\;\;\;
\varepsilon=
\begin{cases}
1 & \text{if $q_1=1$} \\
0 & \text{if $q_1\neq 1$}
\end{cases}.
\end{equation*}
Moreover, we assume that $n_{1+\delta}\leq \cdots\leq n_r$ and
$m_{1+\varepsilon}\leq \cdots\leq m_r$.

Define $b_0\Sigma(\mathbf{n},\mathbf{p})$, $b_1\Sigma(\mathbf{n},\mathbf{p})$ and
$b_0\Sigma(\mathbf{m},\mathbf{q})$, $b_1\Sigma(\mathbf{m},\mathbf{q})$ as following:
\begin{align*}
&b_0\Sigma(\mathbf{n};\mathbf{p})
:=\left\{(\zeta_1,\cdots,\zeta_r)\in\mathbb{C}^{n_1}\times\cdots\times\mathbb{C}^{n_r}:
\sum_{i=1}^r\|\zeta_i\|^{2p_i}=1, \; \|\zeta_j\|\neq 0,1+\delta\leq j\leq r\right \},\\
&b_1\Sigma(\mathbf{n};\mathbf{p}) :=\bigcup_{j=1+\delta}^r
\left\{(\zeta_1,\cdots,\zeta_r)\in\mathbb{C}^{n_1}\times\cdots\times\mathbb{C}^{n_r}:
\sum_{i=1}^r\|\zeta_i\|^{2p_i}=1, \; \|\zeta_j\|= 0\right \},\\
&b_0\Sigma(\mathbf{m};\mathbf{q})
:=\left\{(\xi_1,\cdots,\xi_r)\in\mathbb{C}^{m_1}\times\cdots\times\mathbb{C}^{m_r}:
\sum_{i=1}^r\|\xi_i\|^{2q_i}=1,\; \|\xi_j\|\neq 0,1+\varepsilon\leq
j\leq r\right\},\\
&b_1\Sigma(\mathbf{m};\mathbf{q})
:=\bigcup_{j=1+\varepsilon}^r\left\{(\xi_1,\cdots,\xi_r)\in\mathbb{C}^{m_1}\times\cdots\times\mathbb{C}^{m_r}:
\sum_{i=1}^r\|\xi_i\|^{2q_i}=1,\; \|\xi_j\|= 0\right\}.
\end{align*}
Then $b_0\Sigma(\mathbf{n};\mathbf{p})$(resp.
$b_0\Sigma(\mathbf{m};\mathbf{q})$ ) consists of all strongly
pseudoconvex points of $b\Sigma(\mathbf{n};\mathbf{p})$ (resp.
$b\Sigma(\mathbf{m};\mathbf{q})$ ) and any one of
$b_1\Sigma(\mathbf{n};\mathbf{p})$(resp.
$b_1\Sigma(\mathbf{m};\mathbf{q})$ ) is not a strongly pseudoconvex
point of $b\Sigma(\mathbf{n};\mathbf{p})$ (resp.
$b\Sigma(\mathbf{m};\mathbf{q})$ ). Since $h$ is a biholomorphic
linear isomorphism, we have
\begin{equation}
h(b_0\Sigma(\mathbf{n};\mathbf{p}))=b_0\Sigma(\mathbf{m};\mathbf{q})
\end{equation}
and
\begin{equation}
h(b_1\Sigma(\mathbf{n};\mathbf{p}))=b_1\Sigma(\mathbf{m};\mathbf{q}).
\end{equation}

Let
\begin{eqnarray}
\begin{aligned}
h(\zeta_1,\zeta_2,\cdots,\zeta_r)=(\zeta_1,\zeta_2,\cdots,\zeta_r)
\begin{pmatrix}
D_{11}& D_{12} & \cdots & D_{1r} \\
D_{21}& D_{22} & \cdots & D_{2r} \\
\cdots& \cdots & \cdots & \cdots \\
D_{r1}& D_{r2} & \cdots & D_{rr}
\end{pmatrix}
\end{aligned}.
\end{eqnarray}
According to whether $p_1$ or $q_1$ is equal to $1$ or not, there
are two cases: {\bf (i)} Neither $p_1$ nor $q_1$ equals to $1$; {\bf
(ii)} Either $p_1$ or $q_1$ equals to $1$.

{\bf Case (i).} In this case, we have $p_1\neq 1,\; q_1\neq 1$ and
thus $\delta=0,\; \varepsilon=0$.
\begin{align*}
&b_0\Sigma(\mathbf{n};\mathbf{p})
:=\left\{(\zeta_1,\cdots,\zeta_r)\in\mathbb{C}^{n_1}\times\cdots\times\mathbb{C}^{n_r}:
\sum_{i=1}^r\|\zeta_i\|^{2p_i}=1, \; \|\zeta_j\|\neq 0,1\leq j\leq r \right\},\\
&b_0\Sigma(\mathbf{m};\mathbf{q})
:=\left\{(\xi_1,\cdots,\xi_r)\in\mathbb{C}^{m_1}\times\cdots\times\mathbb{C}^{m_r}:
\sum_{i=1}^r\|\xi_i\|^{2q_i}=1,\; \|\xi_j\|\neq 0,1\leq j\leq r
\right\}.
\end{align*}

Since $h$ a biholomorphic linear isomorphism, we can assume $m_1\leq
n_1$. Hence, we have $m_1\leq n_1\leq \cdots\leq n_r$. In the
following we will use Lemma 2.7 to prove that there exist exactly
one nonzero block in the first column of the matrix of $h$. Consider
the system of linear equations
\begin{equation}
\sum_{j=1}^{r}\zeta_j D_{j1}=0.
\end{equation}
Suppose that $(\alpha_1,\cdots,\alpha_r)$ is a solution of (16) with
$\|\alpha_1\|^2\neq 0,\cdots,\|\alpha_r\|^2\neq 0$. Then there
exists a $\lambda>0$ such that
$$\lambda \sum_{j=1}^{r}\|\alpha_j\|^{2p_j}=1,$$
that is,
$(\lambda^{\frac{1}{2p_1}}\alpha_1,\cdots,\lambda^{\frac{1}{2p_r}}\alpha_r)\in
b_0\Sigma(\mathbf{n};\mathbf{p})$. But the first component of
$h(\lambda^{\frac{1}{2p_1}}\alpha_1,\cdots,\lambda^{\frac{1}{2p_r}}\alpha_r)$
is $\mathbf{0}\in\mathbb{C}^{m_1}$, and thus
$h(\lambda^{\frac{1}{2p_1}}\alpha_1,\cdots,\lambda^{\frac{1}{2p_r}}\alpha_r)\not\in
b_0\Sigma(\mathbf{m};\mathbf{q})$. This is a contradiction with
(13). Thus, the system of linear equations (16) does not have
solution $(\alpha_1,\cdots,\alpha_r)$ with $\|\alpha_1\|^2\neq
0,\cdots,\|\alpha_r\|^2\neq 0$. By Lemma 2.7, there is exactly one
$D_{j_11}\neq \mathbf{0}$. After a permutation $\sigma_1$ of row
index of $D_{ij}$ (which is equivalent to a permutation $\sigma_1$
of the index of $\zeta_i \;(1\leq i\leq r)$ in the (15)), we can
assume that $j_1=1$. Thus, $D_{11}$ is a nonsingular $m_1\times m_1$
matrix with $m_1=n_{\sigma_1(1)}$ and $D_{j1}=\mathbf{0}\;(2\leq
j\leq r)$. Therefore, the first group of components of the mapping
$h$ is independent of the variables
$\zeta_{\sigma_1(2)},\cdots,\zeta_{\sigma_1(r)}$. For the simplicity
of notation, we assume that $\sigma_1$ is the identity permutation,
i.e. $\sigma_1(i)=i\;(1\leq i\leq r)$.

Next, let
\begin{align*}
Pr_1\Sigma(\mathbf{n};\mathbf{p}):=\Sigma(\mathbf{n},\mathbf{p})\cap
\{\zeta_1=0\},\;\;
Pr_1\Sigma(\mathbf{m};\mathbf{q}):=\Sigma(\mathbf{m},\mathbf{q})\cap
\{\xi_1=0\}.
\end{align*}
Since the first group of components of $h$ is independent of
$\zeta_2,\cdots,\zeta_r$, we can consider the restriction $\tilde{h}
$ of $h$ to $\Sigma(\mathbf{n};\mathbf{p})\cap
\{\zeta_1=0\}=:Pr_1\Sigma(\mathbf{n};\mathbf{p})$ as follows:
$$\tilde{h}:Pr_1\Sigma(\mathbf{n};\mathbf{p})\rightarrow Pr_1\Sigma(\mathbf{m};\mathbf{q})$$
\begin{eqnarray*}
\begin{aligned}
\tilde{h}(\zeta_2,\cdots,\zeta_r)=(\zeta_2,\cdots,\zeta_r)\begin{pmatrix}
D_{22} & D_{23} & \cdots & D_{2r} \\
D_{32} & D_{33} & \cdots & D_{3r} \\
\cdots & \cdots & \cdots &\cdots  \\
D_{r2} & D_{r3} & \cdots & D_{rr} \\
\end{pmatrix}.
\end{aligned}
\end{eqnarray*}
Thus $\tilde{h}$ is a biholomorphic linear mapping between
$Pr_1\Sigma(\mathbf{n};\mathbf{p})$and
$Pr_1\Sigma(\mathbf{m};\mathbf{q})$. By the same reasoning as above,
we get that, after a permutation $\sigma_2$ of the index of
$\zeta_i\; (2\leq i\leq r)$, $D_{22}$ is a nonsingular $m_2\times
m_2$ matrix with $m_2=n_{\sigma_2(2)}$ and $D_{j2}=\mathbf{0}$ for
$3\leq j\leq r$. Again, for the simplicity of notation, we assume
that $\sigma_2$ is the identity permutation.

In the same way we can show that for each $i=1,\cdots,r$, after a
permutation $\sigma_i$ of the index of $\zeta_j\;  (i\leq j\leq r)$,
$D_{ii}$ is a nonsingular $m_i\times m_i$ matrix with
$m_i=n_{\sigma_i(i)}$ and $D_{jk}=\mathbf{0}$ for $k<j\leq r$. Thus,
if we let $\sigma=\sigma_r\circ\cdots\circ\sigma_1$, then we have
\begin{eqnarray*}
\begin{aligned}
h(\zeta_1,\zeta_2,\cdots,\zeta_r)
=(\zeta_{\sigma(1)},\zeta_{\sigma(2)},\cdots,\zeta_{\sigma(r)})
\begin{pmatrix}
D_{11}     & D_{12}     & \cdots & D_{1r} \\
\mathbf{0} & D_{22}     & \cdots & D_{2r} \\
\cdots     & \cdots     & \ddots & \cdots \\
\mathbf{0} & \mathbf{0} & \cdots & D_{rr}
\end{pmatrix}
\end{aligned}.
\end{eqnarray*}

Now we prove that $D_{ij}=0\; (1\leq i<j\leq r)$. In fact, we will
show that for $1\leq i\leq r$, the $i$-th column of the above matrix
of $h$ has only one nonzero block. Since every block $D_{ii}\;(1\leq
i\leq r)$ is nonsingular, we get that all other blocks $D_{ij}=0\;
(1\leq i<j\leq r)$.

Suppose that there exist at least two nonzero blocks on some column,
say the last column, of the above matrix of $h$. Then the system of
linear equations
\begin{equation*}
\sum_{j=1}^{r}\zeta_{\sigma(j)}D_{jr}=0
\end{equation*}
has a solution
$(\gamma_{\sigma(1)},\gamma_{\sigma(2)},\cdots,\gamma_{\sigma(r)})\in
(\mathbb{C}^{n_{\sigma(1)}}\setminus\{0\})\times(\mathbb{C}^{n_{\sigma(2)}}\setminus\{0\})\times\cdots\times
(\mathbb{C}^{n_{\sigma(r)}}\setminus \{0\})$ such that
$\sum_{j=1}^r\|\gamma_j\|^{2p_j}=\sum_{j=1}^r\|\gamma_{\sigma(j)}\|^{2p_{\sigma(j)}}=1$.
That is $(\gamma_1,\cdots,\gamma_r)\in
b_0\Sigma(\mathbf{n},\mathbf{p})$, but
$h(\gamma_1,\cdots,\gamma_r)\not\in
b_0\Sigma(\mathbf{m};\mathbf{q})$. This is a contradiction with
(13). Thus, each column of the above matrix of $h$ has only one
nonzero block and $D_{ij}=0,1\leq i<j\leq r$. That is,
\begin{eqnarray*}
\begin{aligned}
h(\zeta_1,\zeta_2,\cdots,\zeta_r)
=(\zeta_{\sigma(1)},\zeta_{\sigma(2)},\cdots,\zeta_{\sigma(r)})
\begin{pmatrix}
 D_{11}     & \mathbf{0} & \cdots & \mathbf{0} \\
 \mathbf{0} & D_{22}     & \cdots & \mathbf{0} \\
 \cdots     & \cdots     & \ddots & \cdots     \\
 \mathbf{0} & \mathbf{0} & \cdots & D_{rr}
\end{pmatrix}
\end{aligned}.
\end{eqnarray*}

For each fixed $j$ ($1\leq j\leq r$), if
$\|\zeta_{\sigma(j)}\|^2<1$, then
$(\underbrace{0,\cdots,0}_{\sigma(j)-1},\zeta_{\sigma(j)},0,\cdots,0)\in
 \Sigma(\mathbf{n},\mathbf{p})$
and $h(\underbrace{0,\cdots,0}_{\sigma(j)-1
},\zeta_{\sigma(j)},0,\cdots,0)
=(\underbrace{0,,\cdots,0}_{j-1},\zeta_{\sigma(j)},0,\cdots,0)
\text{diag}(D_{11},\cdots,D_{rr})
=(0,0,\cdots,0,\zeta_{\sigma(j)}D_{jj},0,$ $\cdots,0) \in
\Sigma(\mathbf{m},\mathbf{q})$. Thus $\|\zeta_{\sigma(j)}
D_{jj}\|^2< 1$. On the other hand, for
$\|\xi_{\sigma^{-1}(j)}\|^2<1$, then
$\|\xi_{\sigma^{-1}(j)}D^{-1}_{jj}\|^2$ $< 1$. This indicates that
$D_{jj}(\mathbf{B}^{n_{\sigma(j)}})\subset
\mathbf{B}^{n_{\sigma(j)}}$ and
$D_{jj}^{-1}(\mathbf{B}^{n_{\sigma(j)}})\subset
\mathbf{B}^{n_{\sigma(j)}}$. Therefore, $D_{jj}$ is a unity
transformation of $\mathbb{C}^{n_{\sigma(j)}}$ for $1\leq j\leq r$.

For $2\leq i\leq r$ and any
$\zeta_{\sigma(1)}\in\mathbb{C}^{n_{\sigma(1)}}$ with
$\|\zeta_{\sigma(1)}\|^2<1$, there exists
$\zeta_{\sigma(i)}\in\mathbb{C}^{n_{\sigma(i)}}$ such that
$(0,\cdots,0,\zeta_{\sigma(1)},0,\cdots,0,\zeta_{\sigma(i)},0,\cdots,0)\in
b\Sigma(\mathbf{n};\mathbf{p})$. Thus,
$(\zeta_{\sigma(1)}D_{11},0,\cdots,0,\zeta_{\sigma(i)}D_{ii},0,\cdots,0)\in
b\Sigma(\mathbf{m};\mathbf{q})$. That is, from
$\|\zeta_{\sigma(1)}\|^{2p_{\sigma(1)}}+\|\zeta_{\sigma(i)}\|^{2p_{\sigma(i)}}=1,$
we can get
$$\|\zeta_{\sigma(1)}\|^{2q_1}+\|\zeta_{\sigma(i)}\|^{2q_i}(=
\|\zeta_{\sigma(1)}D_{11}\|^{2q_1}+\|\zeta_{\sigma(i)}D_{ii}\|^{2q_i})=1.$$
Hence, $p_{\sigma(i)}=q_i,\; 1\leq i\leq r$. This finish the proof
of Lemma 2.7 in {Case (i)}.

{\bf Case (ii).} In this case, without loss of generality (note that
$h:\Sigma(\mathbf{n};\mathbf{p})\rightarrow
\Sigma(\mathbf{m};\mathbf{q})$ is a biholomorphic linear
isomorphism), we can assume $q_1=1$, and then $\varepsilon=1$. We
will prove $p_1=1$ here.

Since $b_1\Sigma(\mathbf{n};\mathbf{p})=\bigcup_{j=1+\delta}^r b
Pr_j\Sigma(\mathbf{n},\mathbf{p})$ and
$b_1\Sigma(\mathbf{m};\mathbf{q})=\bigcup_{j=2}^r b
Pr_j\Sigma(\mathbf{m};\mathbf{q})$, where
$Pr_j\Sigma(\mathbf{n},\mathbf{p}):=\Sigma(\mathbf{n},\mathbf{p})\cap
\{\zeta_j=0\}$ and $Pr_j
\Sigma(\mathbf{m};\mathbf{q}):=\Sigma(\mathbf{m};\mathbf{q})\cap\{\xi_j=0\}$,
by (14), we have
$$h(\bigcup_{j=1+\delta}^r b Pr_j\Sigma(\mathbf{n},\mathbf{p}))\subset
\bigcup_{j=2}^r b Pr_j\Sigma(\mathbf{m};\mathbf{q}).$$ By the same
argument in the proof of Lemma 2.6 above, we can get
\begin{equation*}
h(\mathbf{B}^{n_1}\times\{0\}\times\cdots\times\{0\})\subset
\mathbf{B}^{m_1}\times\{0\}\times\cdots\times\{0\}.
\end{equation*}
Thus we have $D_{1j}=\mathbf{0}\; (2\leq j\leq r)$ for the matrix of
$h$. Apply the same argument to $h^{-1}$, we get
\begin{equation*}
h^{-1}(\mathbf{B}^{m_1}\times\{0\}\times\cdots\times\{0\})\subset
\mathbf{B}^{n_{1}} \times\{0\}\times\cdots\times\{0\}.
\end{equation*}
Thus $h\mid_{\mathbf{B}^{n_1}\times\{0\}\times\cdots\times\{0\}}$ is
a biholomorphism between
$\mathbf{B}^{n_1}\times\{0\}\times\cdots\times\{0\}$ and
$\mathbf{B}^{m_1}\times\{0\}\times\cdots\times\{0\}$. In particular,
we get that $n_1=m_1$.

Since $h$ is a holomorphic linear isomorphism of
$\mathbb{C}^{|\mathbf{n}|}$ onto $\mathbb{C}^{|\mathbf{m}|}$ and
$D_{1j}=0\;(2\leq j\leq r)$, we obtain that $D_{11}$ and
$\begin{pmatrix}
D_{22} & \cdots & D_{2r}\\
\cdots & \cdots & \cdots\\
D_{r2} & \cdots & D_{rr}
\end{pmatrix}$
are invertible constant matrices. Moreover, we have
\begin{eqnarray*}
\begin{aligned}
h^{-1}(\xi_1,\cdots,\xi_r)&=(\xi_1,\cdots,\xi_r)
\begin{pmatrix}D_{11}^{-1} & 0 &\cdots & 0\\
E_{21}& E_{22} &\cdots & E_{2r}\\
\cdots&\cdots &\cdots &\cdots \\
E_{r1}& E_{r2} &\cdots & E_{rr}\\
\end{pmatrix}.
\end{aligned}
\end{eqnarray*}

If $\sum_{j=2}^{r}\|\zeta_j\|^{2p_j}<1$, then
$(0,\zeta_2,\cdots,\zeta_r)\in\Sigma(\mathbf{n},\mathbf{p})$ and
$$h(0,\zeta_2,\cdots,\zeta_r)=(\sum_{j=2}^r\zeta_jD_{j1},\sum_{j=2}^r\zeta_jD_{j2},
\cdots,\sum_{j=2}^r\zeta_jD_{jr})\in
\Sigma(\mathbf{m},\mathbf{q}).$$
Thus
$$\sum_{k=2}^r\|\sum_{j=2}^r\zeta_jD_{jk}\|^{2q_k}
<1-\|\sum_{j=2}^r\zeta_jD_{j1}\|^{2q_1}\leq 1.$$ By the same way,
for $\sum_{j=2}^r\|\xi_j\|^{2q_j}< 1$, we have
$$\sum_{k=2}^r\|\sum_{j=2}^r\xi_jE_{jk}\|^{2p_k}<1-\|\sum_{j=2}^r\xi_jE_{j1}\|^2\leq
1.$$

This indicates that the mapping
$$\tilde{h}:Pr_1\Sigma(\mathbf{n};\mathbf{p})\rightarrow Pr_1\Sigma(\mathbf{m};\mathbf{q})$$
\begin{eqnarray*}
\begin{aligned}
\tilde{h}(\zeta_2,\cdots,\zeta_r)=(\zeta_2,\cdots,\zeta_r)\begin{pmatrix}
D_{22}& \cdots & D_{2r}\\
\cdots&\cdots &\cdots \\
D_{r2} &\cdots & D_{rr}\\
\end{pmatrix}
\end{aligned}
\end{eqnarray*}
is a biholomorphic linear mapping between
$Pr_1\Sigma(\mathbf{n};\mathbf{p})$and
$Pr_1\Sigma(\mathbf{m};\mathbf{q})$. Since $p_2\neq 1$ and $q_2\neq
1$, we can apply the conclusion in the case (i) to get that there
exists a permutation $\sigma \;(\in S_{r-1})$ of $\{2,\cdots, r\}$
such that $n_{\sigma(i)}=m_i,\; p_{\sigma(i)}=q_i \;(2\leq i\leq r)$
and
\begin{eqnarray*}
\begin{aligned}
\tilde{h}(\zeta_2,\cdots,\zeta_r)=(\zeta_{\sigma(2)},\cdots,\zeta_{\sigma(r)})
\begin{pmatrix}
D_{22} &       &    \\
       &\ddots &    \\
       &       & D_{rr}\\
\end{pmatrix},
\end{aligned}
\end{eqnarray*}
where $D_{ii}$ is a unitary transformation of
$\mathbb{C}^{m_i}(m_i=n_{\sigma(i)})$ for $2\leq i\leq r$.
Therefore,
\begin{eqnarray*}
\begin{aligned}
h(\zeta_1,\zeta_2,\cdots,\zeta_r)=(\zeta_1,\zeta_{\sigma(2)},\cdots,\zeta_{\sigma(r)})
\begin{pmatrix}
D_{11} &        &       &       \\
D_{21} & D_{22} &       &       \\
\vdots &        &\ddots &       \\
D_{r1} &        &       & D_{rr}\\
\end{pmatrix}.
\end{aligned}
\end{eqnarray*}

If $\|\zeta_1\|^2<1$, then $(\zeta_1,0,\cdots,0)\in
\Sigma(\mathbf{n},\mathbf{p})$ and
$h(\zeta_1,0,\cdots,0)=(\zeta_1D_{11},0,\cdots, 0)\in
\Sigma(\mathbf{m},\mathbf{q})$. Thus $\|\zeta_1 D_{11}\|^2< 1$. On
the other hand, for $\|\xi_1\|^2<1$, then
$\|\xi_1D^{-1}_{11}\|^2<1$. This indicates that
$D_{11}(\mathbf{B}^{n_1})\subset \mathbf{B}^{n_1}$ and
$D_{11}^{-1}(\mathbf{B}^{n_1})\subset \mathbf{B}^{n_1}$. Therefore,
$D_{11}$ is a unity transformation of $\mathbb{C}^{n_1}$.

If $\|\zeta_{\sigma(j)}\|^2=1$, then
$(0,\cdots,0,\zeta_{\sigma(j)},0,\cdots,0)\in
b\Sigma(\mathbf{n};\mathbf{p}),$ and then
$$h(0,\cdots,0,\zeta_{\sigma(j)},0,\cdots,0)=(\zeta_{\sigma(j)}D_{j1},0,\cdots,
0,\zeta_{\sigma(j)}D_{jj},0,\cdots,0)\in
b\Sigma(\mathbf{m};\mathbf{q}).$$
Thus
$$\|\zeta_{\sigma(j)}
D_{j1}\|^2=1-\|\zeta_{\sigma(j)}D_{jj}\|^{2q_j}=1-\|\zeta_{\sigma(j)}\|^{2q_j}=0$$
(Note $D_{jj}$ is a unitary matrix). Hence $D_{j1}=\mathbf{0}$ for
$2\leq j\leq r$. Thus we get
\begin{eqnarray*}
\begin{aligned}
h(\zeta_1,\zeta_2,\cdots,\zeta_r)=(\zeta_1,\zeta_{\sigma(2)},\cdots,\zeta_{\sigma(r)})
\begin{pmatrix}
D_{11} &        &       &       \\
       & D_{22} &       &       \\
       &        &\ddots &       \\
       &        &       & D_{rr}\\
\end{pmatrix},
\end{aligned}
\end{eqnarray*}
where $D_{ii}$ are unitary transformations of $\mathbb{C}^{n_i}$ for
$1\leq i\leq r$.

Take $\zeta_1\in\mathbb{C}^{n_1},
\zeta_{\sigma(2)}\in\mathbb{C}^{n_{\sigma(2)}}$ (note $\sigma$ is a
permutation of $\{2,\cdots, r\}$) such that
$\|\zeta_1\|^{2p_1}=\frac{1}{2}$ and
$\|\zeta_{\sigma(2)}\|^{2p_{\sigma(2)}}=\frac{1}{2}$. Then
$$h(\underbrace{\zeta_1,0,\cdots,0}_{\sigma(2)-1},
\zeta_{\sigma(2)},0,\cdots,0)=(\zeta_1D_{11},\zeta_{\sigma(2)}D_{22},0,\cdots,0)\in
b\Sigma(\mathbf{m},\mathbf{q}).$$ Hence
$$\|\zeta_1\|^2=\|\zeta_1D_{11}\|^2=1-\|\zeta_{\sigma(2)}D_{22}\|^{2q_2}=
1-\|\zeta_{\sigma(2)}\|^{2q_2}=1-\|\zeta_{\sigma(2)}\|^{2p_{\sigma(2)}}=\frac{1}{2}$$
(Note $q_1=1$ here). So we have $p_1=1\;(=q_1)$. The proof of Lemma
2.8 is finished.

\vskip 8pt Remark. Lemma 2.8 is an extension of Theorem 1.A to the
case of the holomorphic linear isomorphisms between two
equidimensional generalized complex ellipsoids.

\section{Proof of the Main Theorem }

\vskip 8pt
\noindent
{\bf Proof of Theorem 1.1}

\vskip 3pt Let $f:H_{\Omega_1}(\mathbf{n};\mathbf{p})\rightarrow
H_{\Omega_2}(\mathbf{m};\mathbf{q})$ be a biholomorphism. By Lemma
2.6, we have $f(\Omega_1\times\{0\})\subset \Omega_2\times \{0\}$.
In particular, we have $f(0,0)\in \Omega_2\times\{0\}$. Thus, we can
choose an automorphism $\Phi\in
\Gamma(H_{\Omega_2}(\mathbf{m},\mathbf{q}))$ such that $\Phi\circ
f(0,0)=(0,0)$. Thus $\Phi\circ
f:H_{\Omega_1}(\mathbf{n};\mathbf{p})\rightarrow
H_{\Omega_2}(\mathbf{m};\mathbf{q})$ is a biholomorphism with
$\Phi\circ f(0,0)=(0,0)$. Since any Hua domain is a bounded circular
domain and contains the origin, by Cartan's theorem,
$g(\mu,\zeta_1,\cdots,\zeta_r)=\Phi\circ
f(\mu,\zeta_1,\cdots,\zeta_r)$ is a biholomorphic linear mapping
between $H_{\Omega_1}(\mathbf{n};\mathbf{p})$ and
$H_{\Omega_2}(\mathbf{m};\mathbf{q})$, namely
\begin{eqnarray*}
\begin{aligned}
g(z,w_{(1)},\cdots,w_{(r)})=(z,w_{(1)},\cdots,w_{(r)})
\begin{pmatrix}A & B\\
C& D\\
\end{pmatrix}
=(z,w_{(1)},\cdots,w_{(r)})
\begin{pmatrix}A & B_{1} &\cdots & B_{r}\\
C_{1}& D_{11} &\cdots & D_{1r}\\
\cdots&\cdots &\cdots &\cdots \\
C_{r}& D_{r1} &\cdots & D_{rr}\\
\end{pmatrix},
\end{aligned}
\end{eqnarray*}
where $A,B,C,D$ are constant matrices. By Lemma 2.6, we have $g(\Omega_1\times\{0\})=\Omega_2\times \{0\}$,
and this means $B=0$.

Let $|\mathbf{n}|=n_1+\cdots+n_r, |\mathbf{m}|=m_1+\cdots+m_r$. Let
$\Sigma(\mathbf{n};\mathbf{p})$ and $\Sigma(\mathbf{m};\mathbf{q})$
be two generalized complex ellipsoids defined respectively by
\begin{align*}
&\Sigma(\mathbf{n};\mathbf{p}):=\left\{(\zeta_1,\cdots,\zeta_r)\in\mathbb{C}^{n_1}\times\cdots\times\mathbb{C}^{n_r}:
\sum_{j=1}^r\|\zeta_j\|^{2p_j}<1\right \}, \\
&\Sigma(\mathbf{m};\mathbf{q}):=\left\{(\xi_1,\cdots,\xi_r)\in\mathbb{C}^{m_1}\times\cdots\times\mathbb{C}^{m_r}:
\sum_{j=1}^r\|\xi\|^{2q_j}<1\right \}.
\end{align*}

Since $g$ is a holomorphic linear isomorphism of
$\mathbb{C}^{d_1+|\mathbf{n}|}$ onto $\mathbb{C}^{d_2+|\mathbf{m}|}$
and $B=0$, we obtain that $A$ and $D$ are invertible constant
matrices, moreover
\begin{eqnarray*}
\begin{aligned}
g^{-1}(z',w'_{(1)},\cdots,w'_{(r)})&=(z',w'_{(1)},\cdots,w'_{(r)})
\begin{pmatrix}A^{-1} & 0\\
-D^{-1}CA^{-1}& D^{-1}\\
\end{pmatrix}\\
&=(z',w'_{(1)},\cdots,w'_{(r)})
\begin{pmatrix}A^{-1} & 0 &\cdots & 0\\
E_{1}& F_{11} &\cdots & F_{1r}\\
\cdots&\cdots &\cdots &\cdots \\
E_{r}& F_{r1} &\cdots & F_{rr}\\
\end{pmatrix}.
\end{aligned}
\end{eqnarray*}

Note $g(z,w_{(1)},\cdots,w_{(r)})=\Phi\circ
f(z,w_{(1)},\cdots,w_{(r)})$ is a holomorphically linear isomorphism
of $H_{\Omega_1}(\mathbf{n};\mathbf{p})$ to
$H_{\Omega_2}(\mathbf{m};\mathbf{q})$. If
$(\zeta_1,\cdots,\zeta_r)\in \Sigma(\mathbf{n};\mathbf{p}),$ that
is,
$$\sum_{j=1}^{r}\|\zeta_j\|^{2p_j}<1\; (=N_{\Omega_1}(0,0)),$$
then $(0,\zeta_1,\cdots,\zeta_r)\in
H_{\Omega_1}(\mathbf{n},\mathbf{p})$ and
$$g(0,\zeta_1,\cdots,\zeta_r)=(\sum_{j=1}^r\zeta_jC_j,\sum_{j=1}^r\zeta_jD_{j1},\cdots,\sum_{j=1}^r\zeta_jD_{jr})\in
 H_{\Omega_2}(\mathbf{m},\mathbf{q}).$$
Thus, by Proposition 2.1(b), we obtain
$$\sum_{j=1}^r\|\sum_{i=1}^r\zeta_iD_{ij}\|^{2q_j}<
N_{\Omega_2}(\sum_{j=1}^r\zeta_jC_j,
\sum_{j=1}^r\overline{\zeta_jC_j})\leq 1.$$ By the same way, if
$(\xi_1,\cdots,\xi_r)\in\Sigma(\mathbf{m};\mathbf{q})$,then
$$\sum_{j=1}^r\|\sum_{i=1}^r\xi_iF_{ij}\|^{2p_j}<N_{\Omega_1}(\sum_{j=1}^r\xi_jE_j,\sum_{j=1}^r\overline{\xi_jE_j})\leq
1.$$ This indicates that the mapping
$$h:\Sigma(\mathbf{n};\mathbf{p})\rightarrow \Sigma(\mathbf{m};\mathbf{q})$$
\begin{eqnarray*}
\begin{aligned}
h(\zeta_1,\cdots,\zeta_r):=(\zeta_1,\cdots,\zeta_r)\begin{pmatrix}
D_{11}& \cdots & D_{1r}\\
\cdots&\cdots &\cdots \\
D_{r1} &\cdots & D_{rr}\\
\end{pmatrix}
\end{aligned}
\end{eqnarray*}
is a biholomorphic linear mapping between
$\Sigma(\mathbf{n};\mathbf{p})$ and $\Sigma(\mathbf{m};\mathbf{q})$.
By Lemma 2.8, we get
\begin{eqnarray*}
\begin{aligned}
h(\zeta_1,\cdots,\zeta_r)=(\zeta_{\sigma(1)},\cdots,\zeta_{\sigma(r)})
\begin{pmatrix}
D_{11}&       &       \\
      &\ddots &       \\
      &       & D_{rr}\\
\end{pmatrix},
\end{aligned}
\end{eqnarray*}
where $\sigma\in S_r$ is a permutation such that
$n_{\sigma(i)}=m_i$, $p_{\sigma(i)}=q_i$ for $1\leq i\leq r$, and
$D_{ii}$ is a unitary transformation of
$\mathbb{C}^{m_i}(n_{\sigma(i)}=m_i)$ for $1\leq i\leq r$.

Now we prove $C=0$. Thus the matrix of $g=\Phi\circ f$ is a
block diagonal matrix.

Since $g(bH_{\Omega_1}(\mathbf{n};\mathbf{p}))=bH_{\Omega_2}(\mathbf{m};\mathbf{q})$,
we have that if
$$\|\zeta_{\sigma(j)}\|^2=N_{\Omega_1}(0,0)^{\frac{1}{p_{\sigma(j)}}}(=1)$$
(that is, $(\underbrace{ 0,\cdots,
0}_{{\sigma(j)}},\zeta_{\sigma(j)},0,\cdots,0)\in
bH_{\Omega_1}(\mathbf{n},\mathbf{p})$), then
$$g(\underbrace{ 0,\cdots, 0}_{{\sigma(j)}},\zeta_{\sigma(j)},0,\cdots,0)= (\zeta_{\sigma(j)}C_j,0,\cdots,0, \zeta_{\sigma(j)}D_{jj},0,\cdots,0)\in
bH_{\Omega_2}(\mathbf{m},\mathbf{q}).$$ Since $D_{jj}$ is a unitary
matrix, we have
$$\|\zeta_{\sigma(j)}D_{jj}\|^2=N_{\Omega_2}(\zeta_{\sigma(j)}C_j,\overline{\zeta_{\sigma(j)}C_j})^{\frac{1}{q_j}}(=\|\zeta_{\sigma(j)}\|^2=1).$$
By Proposition 2.1(b), we have $\zeta_{\sigma(j)}C_j=0$ for all
$\|\zeta_{\sigma(j)}\|=1$. Thus $C_j=0\; (1\leq j\leq r)$.

Therefore
\begin{eqnarray*}
\begin{aligned}
g(z,w_{(1)},\cdots,\zeta_{(r)})=(z,w_{(\sigma(1))},\cdots,w_{(\sigma(r))})
\begin{pmatrix}A &   &  &  \\
  & D_{11} &  &  \\
  &     &\ddots &  \\
  &     &  & D_{rr}\\
\end{pmatrix},
\end{aligned}
\end{eqnarray*}
where $A$ is a holomorphically linear isomorphism of $\Omega_1$ onto
$\Omega_2$,  $\sigma\in S_r$ is a permutation such that
$n_{\sigma(i)}=m_i$, $p_{\sigma(i)}=q_i$ for $1\leq i\leq r$ and
$D_{ii}$ is a unitary transformation of
$\mathbb{C}^{m_i}(n_{\sigma(i)}=m_i)$ for $1\leq i\leq r$. The proof
of Theorem 1.1 is completed.

\vskip 8pt

\noindent
{\bf Proof of Corollary 1.2}
\vskip 3pt

Since $\Gamma({H_{\Omega}}(\mathbf{n};\mathbf{p}))$ is a subgroup of
$Aut({H_{\Omega}}(\mathbf{n};\mathbf{p})),$ Theorem 1.1 immediately
implies Corollary 1.2. This proves Corollary 1.2.

\vskip 5pt \noindent {\bf Proof of Theorem 1.3} \vskip 3pt Since
$f:H_{\Omega_1}(\mathbf{n}_1;\mathbf{p}_1)\rightarrow
H_{\Omega_2}(\mathbf{n}_2;\mathbf{p}_2)$ is a proper holomorphic
mapping between two equidimensional Hua domains,  by Proposition
2.3, $f$ extends holomorphically to a neighborhood $V$ of
$\overline{H_{\Omega_1}(\mathbf{n}_1;\mathbf{p}_1)}$ with
$$f(bH_{\Omega_1}(\mathbf{n}_1;\mathbf{p}_1))\subset bH_{\Omega_2}(\mathbf{n}_2;\mathbf{p}_2).$$
Define
$$S:=\{\xi\in V:J_f(\xi)=0\},$$
where $J_f(\xi)=\det(\frac{\partial f_i}{\partial \xi_j})(\xi)$ is
the holomorphic Jacobian  determinant of
$$f(\xi)=(f_1(\xi),\cdots,f_{d+|\mathbf{n}|}(\xi))\;(\xi\in V).$$

By Proposition 2.4 and Lemma 2.5, we have
\begin{equation}
f(S\cap
b_0H_{\Omega_1}(\mathbf{n}_1;\mathbf{p}_1))\subset(b\Omega_2\times
\{0\})\cup b_1H_{\Omega_2}(\mathbf{n}_2;\mathbf{p}_2).
\end{equation}
If $S\cap H_{\Omega_1}(\mathbf{n}_1;\mathbf{p}_1)\neq \emptyset$,
then, from the assumption that the subset $(b\Omega_1\times
\{0\})\cup b_1H_{\Omega_1}(\mathbf{n}_1;\mathbf{p}_1)$ of
$bH_{\Omega_1}(\mathbf{n}_1,\mathbf{p}_1)$ is contained in some
complex analytic set of complex codimension at least $2$, we have
$S\cap b_0H_{\Omega_1}(\mathbf{n}_1;\mathbf{p}_1)\neq\emptyset$.
Take an irreducible component $S'$ of $S$ with $S'\cap
H_{\Omega_1}(\mathbf{n}_1;\mathbf{p}_1)\neq \emptyset$. Then, from
Proposition 2.4 (a),  the intersection $E_{S'}$ of $S'$ with
$b_0H_{\Omega_1}(\mathbf{n}_1;\mathbf{p}_1)$ is a real analytic
submanifold of real dimension $2(d_1+|\mathbf{n}_1|)-3$ on a dense,
open subset of $E_{S'}$ (Otherwise, $S'$ cannot be separated by
$E_{S'}$ (e.g., see Rudin \cite{R}, Theorem 14.4.5) and thus $S'$
cannot be separated by $bH_{\Omega_1}(\mathbf{n}_1;\mathbf{p}_1).$
This is impossible). From (17), we also have $$f(E_{S'})\subset
(b\Omega_2\times \{0\})\cup
b_1H_{\Omega_2}(\mathbf{n}_2;\mathbf{p}_2).$$
 Thus, by the
uniqueness theorem,
\begin{equation}
f(S'\cap
H_{\Omega_1}(\mathbf{n}_1;\mathbf{p}_1))\subset\bigcup_{j=1+\varepsilon}^r
Pr_j (H_{\Omega_2}(\mathbf{n}_2;\mathbf{p}_2)),
\end{equation}
where
$Pr_j(H_{\Omega_2}(\mathbf{n}_2;\mathbf{p}_2)):=H_{\Omega_2}(\mathbf{n}_2;\mathbf{p}_2)\cap
\{w'_{(j)}=0 \}$ for $1+\varepsilon\leq j\leq r$ and
\begin{equation*}
\varepsilon=
\begin{cases}
1 & \text{if the first component $p_2^1$ of $\mathbf{p}_2$ equals to 1}, \\
0 & \text{if the first component $p_2^1$ of $\mathbf{p}_2$ does not
equal to 1}.
\end{cases}
\end{equation*}
Since codim$S'=1$, codim$(\bigcup_{j=1+\varepsilon}^r Pr_j
(H_{\Omega_2}(\mathbf{n}_2;\mathbf{p}_2))) =
\min\{n_2^{1+\varepsilon},\cdots,n_2^r\}\geq 2$, where
$\mathbf{n}_2=(n_2^1,\cdots,n_2^r)$. Since
$f:H_{\Omega_1}(\mathbf{n}_1;\mathbf{p}_1)\rightarrow
H_{\Omega_2}(\mathbf{n}_2;\mathbf{p}_2)$ is proper, this is a
contradiction with (18). This means $S\cap
H_{\Omega_1}(\mathbf{n}_1;\mathbf{p}_1)=\emptyset$.

Thus $f:H_{\Omega_1}(\mathbf{n}_1;\mathbf{p}_1)\rightarrow
H_{\Omega_2}(\mathbf{n}_2;\mathbf{p}_2)$ is unbranched. Since Hua
domain is simply connected by Proposition 2.1(c), we get that
$f:H_{\Omega_1}(\mathbf{n}_1;\mathbf{p}_1)\rightarrow
H_{\Omega_2}(\mathbf{n}_2;\mathbf{p}_2)$ is a biholomorphism. The
proof of Theorem 1.3 is completed.

\vskip 5pt \noindent {\bf Proof of Corollary 1.4} \vskip 3pt From
Theorem 1.3 and Corollary 1.2, we immediately get Corollary 1.4.

\vskip 5pt \noindent {\bf Proof of Corollary 1.5} \vskip 3pt From
Theorem 1.3, Corollary 1.5 is obviously.

\vskip 6pt

\noindent\textbf{Acknowledgments}\quad The authors are grateful to
Professors Ngaiming Mok and Xiaojun Huang for helpful suggestions.
In addition, the authors thank the referees for many useful
comments. The project is supported by the National Natural Science
Foundation of China (No.11271291).


\begin{thebibliography}{40}
\bibitem{ABP}H. Ahn,  J. Byun,  J. Park, Automorphisms of the Hartogs type
domains over classical symmetric domains, Internat. J. Math. {\bf
23}(No.9)(2012), 1250098, 11 pp.

\bibitem{Ale}H. Alexander, Proper holomorphic mappings in $\mathbb{C}^n,$ Indiana Univ. Math. J. {\bf 26}(1977), 137-146.

\bibitem{Bedf}E. Bedford and S. Bell, Proper self-maps of weakly pseudoconvex domains, Math. Ann. {\bf 261}(1982), 47-49.

\bibitem{Bell}S. Bell, Proper holomorphic mappings between circular domains, Comment. Math. Helvetici {\bf 57}(1982), 532-538.

\bibitem{Bell82}S. Bell, The Bergman kernel function and proper holomorphic mappings, Trans. Amer. Math. Soc. 270(1982), 685-691.

\bibitem{Dini91}G. Dini and A. Primicerio, Proper holomorphic mappings between generalized pseudoellipsoids,
Annali di matematica pura ed applicata {\bf 158}(1991), 219-229.

\bibitem{Dini}G. Dini and A. Primicerio, Localization principle of automorphisms on generalized pseudoellipsoids,
J. Geom. Anal. {\bf 7}(4)(1997), 575-584.


\bibitem{Tu3}Z.M. Feng and Z.H. Tu, On canonical metrics on
Cartan-Hartogs domains, Math. Z. 278(2014), 301-320.


\bibitem{Forst}F. Forstneric, Proper holomorphic mappings: A survey, in Several Complex Variables (edited by J. Fornaess),
Math Notes Vol. {\bf 38}, Princeton University Press, 1993, 297-363.


\bibitem{H} H. Hamada, On proper holomorphic self-maps of
generalized complex ellipsoids, J. Geom. Anal. {\bf 8}(1998),
441-446.


\bibitem{Henkin}G.M. Henkin and R. Novikov, Proper mappings of classical domains, in Linear and Complex Analysis
Problem Book, Lecture Notes in Math. Vol. {\bf 1043}, Springer,
Berlin, 1984, 625-627.

\bibitem{Hua}L.K. Hua, Harmonic Analysis of Functions of Several Complex Variables in the Classical Domains,
Amer. Math. Soc., Providence, RI, 1963.


\bibitem{H} X.J. Huang, On a linearity problem for proper holomorphic maps
between balls in complex spaces of different dimansions, J. Diff.
Geom. {\bf 51}(1999), 13-33.


\bibitem{K1} A. Kodama, On the holomorphic automorphism group of a generalized
complex ellipsoid, Complex Variables and Elliptic Equations {\bf
59}(2014), 1342-1349.

\bibitem{Kodama}A. Kodama, S.G. Krantz and D. Ma, A characterization of generalized complex ellipsoids in $\mathbb{C}^N$
and related results, Indiana Univ. Math. J. {\bf 41}(1992), 173-195.


\bibitem{L} M. Landucci, On the proper holomorphic equivalence for a class of
pseudoconvex domains, Trans. Amer. Math. Soc. {\bf 282}(1984),
807-811.


\bibitem{LZ}A. Loi and M. Zedda, Balanced metrics on Cartan and Cartan-Hartogs domains. Math. Z. \textbf{270}(2012),
1077-1087.

\bibitem{Mok}N. Mok, Extension of germs of holomorphic isometries up to normalizing constants with respect to the Bergman metric,
J. Eur. Math. Soc. {\bf 14}(2012), 1617-1656.

\bibitem{MNT} N. Mok, S.C. Ng and Z.H. Tu, Factorization of proper holomorphic maps on
irreducible bounded symmetric domains of rank $\geq  2$, Sci. China
Math. {\bf 53}(3)(2010), 813-826.

\bibitem{Mok-Tsai}N. Mok and I.H. Tsai, Rigidity of convex realizations of irreducible bounded symmetric domains of rank $\geq 2$,
J. Reine. Angew. Math. {\bf 431}(1992), 91-122.

\bibitem{Naruki}I. Naruki, The holomorphic equivalence problem for a class of Reinhardt domains,
Publ. Res. Inst. Math. Sci., Kyoto Univ. {\bf 4}(1968), 527-543.

\bibitem{Pin}S.I. Pin\v{c}uk, On the analytic continuation of biholomorphic mappings, Math. USSR Sb. {\bf 27}(1975), 375-392.


\bibitem{R}F. Rong, On automorphism groups of generalized Hua domains, Math. Proc. Camb. Phil. Soc. {\bf 156}(2014), 461-472.

\bibitem{R1}W. Rudin, Function Theory in the Unit Ball of
$\mathbb{C}^n$,  Springer-Verlag, 1980.

\bibitem{Tsai}I-H. Tsai, Rigidity of proper holomorphic maps between symmetric domains, J. Diff. Geom. {\bf 37}(1993), 123-160.

\bibitem{Tu1}Z.H. Tu, Rigidity of proper holomorphic maps between equidimensional bounded symmetric domains,
Proc. Amer. Math. Soc. {\bf 130}(2002), 1035-1042.

\bibitem{Tu2}Z.H. Tu, Rigidity of proper holomorphic mappings between nonequidimensional bounded symmetric domains,
Math. Z. {\bf 240}(2002), 13-35.

\bibitem{T-W}Z.H. Tu and L. Wang, Rigidity of proper holomorphic mappings between
certain unbounded non-hyperbolic domains, J. Math. Anal. Appl. {\bf
419}(2014), 703-714.


\bibitem{Tuma}A.E. Tumanov and G.M. Henkin, Local characterization of holomorphic automorphisms of classical domains,
Dokl. Akad. Nauk SSSR {\bf 267}(1982), 796-799. (in Russian)

\bibitem{RWYZ}A. Wang, W.P. Yin, L.Y. Zhang  and G.  Roos, The K\"{a}hler-Einstein metric for some Hartogs domains over bounded symmetric domains,
Science in China Series A: Math. \textbf{49}(9)(2006), 1175-1210.

\bibitem{Yin}W.P. Yin, The Bergman Kernels on Cartan-Hartogs domains, Chinese Sci. Bull. \textbf{44}(21)(1999),
1947-1951.

\bibitem{Yin1}W. Yin, A. Wang, Z. Zhao, X. Zhao and B. Guan, The Bergman kernel functions on Hua domains,
Science in China Series A: Math. {\bf 44}(6)(2001), 727-741.

\end{thebibliography}
\end{document}